\documentclass[11pt]{amsart}

\usepackage[usenames]{color}

\usepackage[colorlinks=true, linkcolor=blue, anchorcolor=blue, citecolor=blue, filecolor=blue, menucolor= blue, urlcolor=red]{hyperref}

\usepackage{amsmath,amssymb,amsthm}

\usepackage{graphicx}
\usepackage[all]{xy}\CompileMatrices
\usepackage{vmargin}
\setpapersize{USletter}
\setmargrb{1in}{1in}{1in}{1in} 

\newtheorem{thm}{Theorem}[subsection]
\newtheorem{prop}[thm]{Proposition}
\newtheorem{lem}[thm]{Lemma}
\newtheorem{cor}[thm]{Corollary}

\theoremstyle{definition}

\newtheorem{defn}[thm]{Definition}
\newtheorem*{notation*}{Notation}
\newtheorem{notation}[thm]{Notation}
\newtheorem{example}[thm]{Example}
\newtheorem{remark}[thm]{Remark}

\newtheorem{rem}[thm]{Remark}

\newtheorem*{remark*}{Remark}

\newcommand{\GMS}{GMS}

\newcommand{\PPP}{{\bf P}}

\newcommand{\pr}[1]{{{\bf P}^{#1}}}
\newcommand{\field}{K}

\newcommand{\defining}[1]{\textbf{#1}}
\DeclareMathOperator{\diag}{diag}
\newcommand{\diff}{\Delta}

\DeclareMathOperator{\reg}{reg}

\DeclareMathOperator{\Image}{Im}
\DeclareMathOperator{\cok}{cok}
\DeclareMathOperator{\codim}{codim}

\newcommand{\card}[1]{\#{#1}}

\renewcommand{\O}{\mathcal{O}}
\newcommand{\I}{\mathcal{I}}
\newcommand{\R}{\mathcal{R}}
\newcommand{\Q}{\mathcal{Q}}

\newcommand{\dvec}{{\bf d}}
\newcommand{\vvec}{{\bf v}}
\newcommand{\avec}{{\bf a}}
\newcommand{\mvec}{{\bf m}}

\begin{document}

\ \vskip-1in
\title[Bounding Hilbert functions of fat points in projective space]{Combinatorial bounds on Hilbert functions of fat points in projective space}

\author[S. Cooper]{Susan Cooper}
\address{Mathematics Department\\
California Polytechnic State University\\
San Luis Obispo, CA 93407, USA \\
(currently Department of Mathematics\\
University of Nebraska\\
Lincoln, NE 68588-0130 USA)
}
\email{sucooper@calpoly.edu, scooper4@math.unl.edu}

\author[B. Harbourne]{Brian Harbourne}
\address{Department of Mathematics\\
University of Nebraska\\
Lincoln, NE 68588-0130 USA}
\email{bharbour@math.unl.edu}

\author[Z. Teitler]{Zach Teitler}
\email{zteitler@boisestate.edu}
\address{Department of Mathematics\\
Boise State University\\
1910 University Drive\\
Boise, ID 83725-1555 USA}

\date{December 4, 2010}

\thanks{Acknowledgments: We thank the NSF-AWM Mentoring Travel Grant for
its support of Cooper, NSA for its support of
Harbourne, and Nebraska IMMERSE
for its support of Teitler. We also thank A. V. Geramita for helpful discussions
which led to the results on graded Betti numbers, J.\ Migliore for several 
illuminating consultations and the referee for his helpful comments.}



\begin{abstract}
We study Hilbert functions of certain non-reduced schemes $A$
supported at finite sets of points 
in $\pr N$, in particular, fat point schemes.
We give combinatorially defined upper and lower bounds for the Hilbert function of $A$
using nothing more than the multiplicities of the points and
information about which subsets of the points are linearly dependent.
When $N=2$, we give these bounds explicitly and we give a sufficient criterion 
for the upper and lower bounds to be equal.
When this criterion is satisfied, we give both a simple
formula for the Hilbert function and
combinatorially defined upper and lower bounds on the graded Betti numbers
for the ideal $I_A$ defining $A$, generalizing results
of Geramita--Migliore--Sabourin~\cite{refGMS}.
We obtain the exact Hilbert functions and graded Betti numbers
for many families of examples, interesting combinatorially, geometrically, and algebraically.
Our method works in any characteristic.
\end{abstract}

\maketitle

\renewcommand{\thethm}{\thesection.\arabic{thm}}
\setcounter{thm}{0}

\section{Introduction}\label{section: intro}

The determination of Hilbert functions and graded Betti numbers for ideals of fat point schemes
in $\pr N$ are issues of significant interest in algebraic geometry, even when $N=2$ (and related areas; e.g., 
\cite{refAP, refHS, refMi}).
Previous work has focused on obtaining exact determinations of Hilbert functions 
and graded Betti numbers in various circumstances. For example,
\cite{refVanc, refGi, refHi2, refHHF} formulate conjectures for the ideals of fat points in $\pr2$.
Various special cases have attracted a lot of attention, such as
when the number of points is a square \cite{refCM3, refEv, refEv2, refHR, refR},
or when the multiplicities of the points are small \cite{refCM1, refGMS, refI} 
or the number of points is small \cite{refNag, refFHH, refGHM}. 

Our aim here is to give a {\it general\/} iterative method for determining  
{\it upper and lower bounds\/} on the Hilbert function of any fat point scheme in $\pr N$
(in any characteristic), using only the multiplicities of the points and information about which subsets of
the points are linearly dependent (meaning, when $N=2$, which subsets of the points are collinear);
i.e., information which can be recovered from the matroid of linear dependencies of the points.
The bounds we obtain in case $N=2$ are explicit and easy to compute.
Although our bounds are tightest when the points are fairly special,
the best results for general sets of points depend on understanding specializations
(see, for example, \cite{refCM2}, \cite{refDu}, \cite{refHi}), so we expect our work to be useful quite broadly.
When $N=2$, we also give a simple criterion for the bounds to coincide, in which case our results determine the
Hilbert function uniquely.
This condition is satisfied in many cases of interest; see \textsection \ref{examplessection} for examples.
When the condition is satisfied, we also give upper and lower bounds on the graded Betti numbers
of the ideal of the fat point scheme, and we give a condition for these bounds to coincide.

We begin with the very general statement Theorem~\ref{genthm},
then specialize to hyperplanes in $\pr N$,
then specialize further to lines in $\pr 2$.

It is an open question to determine the Hilbert functions for specific kinds of
fat point schemes in projective space.
For example, it is not known which functions occur as Hilbert functions
of double point subschemes of $\pr2$. (A double point scheme is a
fat point scheme consisting of a finite set of points all taken with multiplicity 2; i.e., the scheme 
defined by the symbolic square of the ideal defining the reduced finite set of points.)
A start on this problem for points, in the case of double point subschemes of $\pr2$ in 
characteristic $0$, was made in \cite{refGMS} using linkage methods.
One of the main results of \cite{refGMS} 
gives a criterion characterizing a class of functions
each of which occurs as the Hilbert function of a double point scheme.
For each function $h$ in that class, the approach of \cite{refGMS} is to use
a sequence of basic double links to construct a double point scheme 
whose Hilbert function is $h$.

We take an approach opposite to that presented in \cite{refGMS}. For the case of fat point subschemes of $\pr2$,
our procedure amounts to tearing down the subscheme as a sequence of residuals with respect to lines
containing various subsets of points of the scheme to obtain upper and lower bounds
for the Hilbert function of the scheme
(for subschemes $A,B\subseteq \pr N$, by the residual of $A$ with respect to $B$
we mean the subscheme defined by the ideal $I_A:I_B$). 
Starting with any fat point subscheme $Z=Z_0$, 
we choose a sequence of lines
$L_1,\ldots,L_r$ and define $Z_i$ to be the residual of $Z_{i-1}$ with respect
to $L_i$ and we define the associated reduction vector
$\dvec=(d_1,\ldots,d_r)$ by taking $d_i=\deg(L_i\cap Z_{i-1})$.
Given $\dvec$ (or indeed any vector $\dvec = (d_1,\dots,d_n)$ with non-negative integer entries),
we also define functions $f_\dvec$ and $F_\dvec$ 
from the non-negative integers to the non-negative integers
(see Definition \ref{thefunctions1}); these functions
are defined completely in terms of $\dvec$.
Our procedure works in a similar way for fat point subschemes of $\pr N$ for $N>2$, except the 
residuals are taken with respect to hyperplanes, and our bounds are less explicit, since 
the procedure expresses the bounds for a subscheme $Z$ of $\pr N$ in terms of bounds obtained iteratively
on the Hilbert functions of $Z$ and its residuals intersected with various hyperplanes.
The iteration terminates with $\pr 2$ by applying Theorem \ref{mainthmintro}, 
which is proved in section \textsection\ref{bettisection}, ultimately as a corollary of Theorem \ref{genthm}. 
AWK scripts implementing the following theorem
can be obtained at \url{http://www.math.unl.edu/~bharbourne1/CHT/Example.html}.
(Note: in the original posted version \cite{refCHT} of this paper 
the indexation of the entries of $\dvec$ is the reverse of what we use here.) 

\begin{thm}\label{mainthmintro}
Let $Z=Z_0$ be a fat point scheme in $\pr2$
with reduction vector $\dvec = (d_1,\dots,d_r)$ such that $Z_{r+1}=\varnothing$.
Then the Hilbert function $h_Z$ is bounded by $f_\dvec \leq h_Z \leq F_\dvec$.
Furthermore, if $\dvec$ is non-increasing and positive, i.e., $d_1 \geq d_2 \geq \cdots \geq d_r> 0$,
then $f_\dvec = F_\dvec$ unless $\dvec$
contains a subsequence of consecutive entries of the form $(a,a,a)$, or of the form $(a_i,\ldots,a_{i+j+1})$
for $j>1$ where $a_i=a_{i+1}$, $a_{i+j}=a_{i+j+1}$, and $a_{i+1},\ldots,a_{i+j}$ are consecutive integers.
If $d_1 > d_2 > \cdots > d_r> 0$, then $\dvec$ also uniquely determines the graded Betti numbers
of the ideal $I_Z$.
\end{thm}

We will leave an analysis of the bounds obtained by our iterative procedure when $N>2$
to future articles. For the case of $N=2$ our examples in \textsection\ref{examplessection} 
show that this theorem gives the Hilbert functions and graded Betti numbers
of fat point subschemes of $\pr2$ in many situations where they were not previously known,
in some cases for subschemes defined by arbitrary symbolic powers of ideals of points,
substantially extending some of the results of \cite{refGMS} which considered only symbolic squares.

Although the results of \cite{refGMS} were our main motivation, the methods of our paper were inspired
by those of \cite{refFL} and \cite{refFHL}.
In these papers the authors focus their attention on fat point subschemes of $X\subset\pr N$
whose reduced subscheme $X_{\text{red}}$ is contained in a hyperplane $H$.
In this setting, the basic idea of \cite{refFL} and \cite{refFHL} is to obtain the minimal free resolution of the ideal 
$I_X$ (and hence the Hilbert function and
graded Betti numbers for $I_X$) in terms of the Hilbert functions and graded Betti numbers
of $X_i\cap H$, where $\varnothing=X_r\subset X_{r-1}\subset \cdots \subset X_0=X$ are 
successive residuals of $X$ with respect to $H$.

Not only can our method be used to obtain bounds for Hilbert functions of fat point schemes in $\pr N$ for $N > 2$,
in terms of residuation with respect to hyperplanes, but our underlying method can also be applied taking residuals
with respect to nonlinear subvarieties, such as curves in $\pr 2$ of degrees greater than 1,
or with curves (especially smooth rational curves) on surfaces $X$ obtained by blowing up points of $\pr2$. 
We focus on the case of lines in $\pr2$, partly because the
data structures for specifying linear dependencies of points in $\pr2$ are easier to analyze than those for points in higher dimensions,
but also because the work of Geramita, Migliore and Sabourin \cite{refGMS} has sparked interest in finding conditions which 
uniquely determine the Hilbert functions of fat point subschemes of $\pr2$
(as we do in Theorem \ref{GMSthm}, Proposition \ref{GMSprop} and Remark \ref{thefunctions2}).

\renewcommand{\thethm}{\thesubsection.\arabic{thm}}
\setcounter{thm}{0}

\subsection{The key idea}\label{subsection key}
Our fundamental idea is very simple but it ends up giving remarkably strong results.
Given a sheaf $\Q_0$ on a scheme $S$, one often attempts to obtain information about
$h^0(S, \Q_0)$ by finding a short exact sequence 
\addtocounter{thm}{1}
\begin{equation}\label{fundses}
0\to \Q_1\to \Q_0\to \R_0\to 0
\end{equation}
of sheaves where one knows something about the cohomology of $\Q_1$ and $\R_0$.
In particular, the short exact sequence of sheaves gives a long exact sequence of cohomology 
from which we immediately obtain the following bounds:
\addtocounter{thm}{1}
\begin{equation}\label{1stepbnds}
\max\big(h^0(S, \Q_1),h^0(S, \Q_1)+h^0(S, \R_0)-h^1(S, \Q_1)\big)\leq h^0(S, \Q_0)\leq h^0(S, \Q_1)+h^0(S, \R_0).
\end{equation}
In our situation we will know everything about $\R_0$ and we will know 
all we need to know about $h^0(S, \Q_1)-h^1(S, \Q_1)$, but we will not usually know $h^0(S, \Q_1)$ by itself.
However, we will be able to find another short exact sequence in which $\Q_1$ becomes the middle term.
By iteration, we obtain the following theorem (in our typical applications, we will know the quantities $h^0(\Q_i)-h^1(\Q_i)$,
we will either know the quantities $h^0(\R_i)$ or have inductively obtained bounds on them,
and $\Q_{n+1}$ will be a sheaf whose cohomology we know):

\begin{thm}\label{genthm}
Given short exact sequences of sheaves on a scheme $S$
\addtocounter{thm}{1}
\begin{equation}\label{fundexactseqs}
\begin{array}{@{0 \hspace{\tabcolsep} \to \hspace{\tabcolsep} } c @{\hspace{\tabcolsep} \to \hspace{\tabcolsep}} c @{\hspace{\tabcolsep} \to \hspace{\tabcolsep}} c @{\hspace{\tabcolsep} \to \hspace{\tabcolsep}} c }
 \Q_{1} & \Q_0 & \R_0 & 0\\
 \Q_{2} & \Q_{1} & \R_{1} & 0\\
 \multicolumn{3}{c}{\vdots} \\
 \Q_{n+1} & \Q_n & \R_n & 0,
\end{array}
\end{equation}
we obtain the following bounds (in which we have suppressed $S$ in our notation):
\begin{equation*}
\begin{split}
h^0(\Q_0) & \leq h^0(\Q_{n+1}) + h^0(\R_n)+\cdots+ h^0(\R_0)
\end{split}
\end{equation*}
and 
\begin{equation*}
\begin{split}
h^0(\Q_0) & \geq \max\big(h^0(\Q_{n+1}),h^0(\R_n)+h^0(\Q_{n+1})-h^1(\Q_{n+1}),\ldots,h^0(\R_0)+h^0(\Q_1)-h^1(\Q_1)\big).
\end{split}
\end{equation*}
Moreover, the upper bound is an equality if and only if each of the short exact sequences of \eqref{fundexactseqs}
is exact on global sections, and both the upper and lower bounds are equalities 
if $h^0(\R_i)h^1(\Q_{i+1})=0$ for $0\leq i\leq n$.
\end{thm}

\begin{proof}
By applying the bounds in \eqref{1stepbnds} to successive rows of \eqref{fundexactseqs},
we get the upper bound 
\addtocounter{thm}{1}
\begin{equation}\label{eqn: ideal upper bounds}
\begin{split}
h^0(\Q_0) & \leq h^0(\Q_1) + h^0(\R_0) \\
& \leq h^0(\Q_2) + h^0(\R_1)+ h^0(\R_0) \\
& \leq \cdots \\
& \leq h^0(\Q_{n+1}) + h^0(\R_n)+\cdots+ h^0(\R_0)
\end{split}
\end{equation}
and similarly we get the lower bound
\addtocounter{thm}{1}
\begin{equation}\label{eqn: ideal lower bounds}
\begin{split}
h^0(\Q_0) & \geq \max\big(h^0(\Q_1),h^0(\R_0)+h^0(\Q_1)-h^1(\Q_1)\big) \\
& \geq \max\big(h^0(\Q_2),h^0(\R_1)+h^0(\Q_2)-h^1(\Q_2),h^0(\R_0)+h^0(\Q_1)-h^1(\Q_1)\big) \\
& \geq \cdots \\
& \geq \max\big(h^0(\Q_{n+1}) ,h^0(\R_n)+h^0(\Q_{n+1})-h^1(\Q_{n+1}),\ldots,h^0(\R_0)+h^0(\Q_1)-h^1(\Q_1)\big).
\end{split}
\end{equation}
Clearly, the upper bound in \eqref{1stepbnds} is an equality if and only if \eqref{fundses}
is exact on global sections. It follows that the inequalities in \eqref{eqn: ideal upper bounds} are
all equalities if and only if each of the short exact sequences of \eqref{fundexactseqs}
is exact on global sections. In particular, if $h^0(\R_i)h^1(\Q_{i+1})=0$ for
every $0\le i\le n$, then \eqref{eqn: ideal upper bounds} is an equality.
Similarly, $h^0(\Q_i) =\max(h^0(\Q_{i+1}),h^0(\R_i)+h^0(\Q_{i+1})-h^1(\Q_{i+1}))$ if $h^0(\R_i)h^1(\Q_{i+1})=0$, and hence
the inequalities of \eqref{eqn: ideal lower bounds} are all equalities if $h^0(\R_i)h^1(\Q_{i+1})=0$
for $0\leq i\leq n$.
\end{proof}

\begin{remark}\label{equalityrem}
If by induction 
we have upper and lower bounds on $h^0(\R_i)$ and if $h^0(\Q_{n+1})$ and the differences 
$h^0(\Q_i)-h^1(\Q_i)$ are known, then as a corollary of Theorem \ref{genthm} 
we obtain upper and lower bounds on
$h^0(\Q_0)$. This is how we obtain upper and lower bounds on the Hilbert function of fat point subschemes of
projective space in general, and it is how we obtain $f_\dvec$ and $F_\dvec$ in Theorem \ref{mainthmintro} in particular.
Moreover, although we will not usually know $h^1(\Q_{i+1})$, if 
we have $h^0(\R_i)h^1(\Q_{i+1})=0$ for all $i$, then we have
\[
\begin{split}
\max\big(&h^0(\Q_{n+1}) ,h^0(\R_n)+h^0(\Q_{n+1})-h^1(\Q_{n+1}),\ldots,h^0(\R_0)+h^0(\Q_1)-h^1(\Q_1)\big)=\\
& h^0(\Q_0)=h^0(\Q_{n+1}) + h^0(\R_n)+\cdots+ h^0(\R_0)
\end{split}
\]
and thus we do know $h^0(\Q_0)$. 
In particular, if $h^1(\Q_{n+1})=0$ and if for some $i$ we have $h^1(\R_j)=0$ for $j>i$, then it follows from
\eqref{fundexactseqs}, working our way up from the bottom of the diagram, that 
$h^1(\Q_j)=0$ for $j>i$. Hence
if whenever $h^0(\R_i)>0$ we have $h^1(\R_j)=0$ for $j>i$, then 
we have $h^0(\R_i)h^1(\Q_{i+1})=0$ for all $i$, in which case
our upper and lower bounds are both equal to $h^0(\Q_0)$.
This is precisely what happens in those cases where Theorem \ref{mainthmintro} 
allows us to conclude $f_\dvec=F_\dvec$.
\end{remark}

\subsection{Background}\label{Background subsection}
We briefly cover the background we will need regarding
fat points, Hilbert functions, graded Betti numbers and residuals.
For this purpose, let $R=\field[\pr N]$ be the homogeneous coordinate ring for
$\pr N$ over an arbitrary algebraically closed field $\field$.

\begin{defn}\label{fatptdef}
Given points $p_1,\dots,p_r\in\pr N$ and non-negative integers $m_i$,
the ideal $I=\bigcap_{i}I_{p_i}^{m_i}\subset R$, where $I_{p_i}$ is the ideal
generated by all forms that vanish at $p_i$, defines a subscheme $Z\subset \pr N$ 
referred to as a \defining{fat point scheme}. 
\end{defn}

We denote $Z$ by $m_1p_1+\dots+m_rp_r$ and $I$ by $I_Z$.
The \defining{support} of $Z$ is the set of points $p_i$ with $m_i>0$.
If $X=p_1+\dots+p_r$ and $Z=mp_1+\dots+mp_r$, we write $Z=mX$. 

The ideal $I_Z$ is graded in the usual way by degree. 
The span of the forms in $I_Z$ of degree $t$ is denoted $(I_Z)_t$;
thus $I_Z=\bigoplus_t (I_Z)_t$. 
With $(R/I_Z)_t=R_t/(I_Z)_t$, we see that the quotient $R/I_Z$ is also graded.
The \defining{Hilbert function $h_{I_Z}$ of $I_Z$} gives the $\field$-vector space dimension
of $(I_Z)_t$ as a function of $t\ge0$; i.e., $h_{I_Z}(t)=\dim_\field((I_Z)_t)$.
Similarly, the \defining{Hilbert function $h_{Z}$ of $Z$} gives the $\field$-vector space dimension
of $(R/I_Z)_t$ as a function of $t\ge0$; i.e., $h_{Z}(t)=\dim_\field((R/I_Z)_t)$.
It is convenient to regard $h_Z$ as the sequence $(h_Z(0), h_Z(1), \dots)$ and vice versa.

As a graded ideal, $I_Z$ has a minimal graded free resolution. 
This is an exact sequence of the form
$0\to F_1\to F_0\to I_Z\to 0$, where each $F_i$ is a graded free $R$-module and the map
$F_1\to F_0$ is given by a matrix of homogeneous forms of positive degree.
As graded free modules, each $F_i$ is, up to graded isomorphism, 
either $0$ or a direct sum of modules of the form $R[-j]$, where
$R[-j]$ is just $R$ with a shift in the grading; i.e., 
$R[-j]_t=R_{t-j}$.
The \defining{graded Betti numbers} of $I_Z$ are the numbers of factors of $R[-j]$ in each sum.
In particular, we can write $F_0$ as $\bigoplus_{j\geq0} R[-j]^{\nu_j(Z)}$ and
$F_1$ as $\bigoplus_{j\geq0} R[-j]^{\sigma_j(Z)}$. 
The graded Betti numbers of $I_Z$ 
are just the sequences $\nu(Z)$ and $\sigma(Z)$.
(For clarity we suppress the $Z$ if it is not needed.)
The value $\nu_j$ specifies the number of 
generators of $I_Z$ in degree $j$ in any minimal set of homogeneous generators of $I_Z$;
the sequence $\sigma$ does the same for the syzygies. In cases in which the Hilbert function
$h_Z$ is known, one can recover the sequence $\sigma$ from $h_Z$ and $\nu$ using the formula
(see \cite{refFHH}) $\sigma_t-\nu_t=-\Delta^3h_{I_Z}(t)$
(which for $t>0$ can be written as $\sigma_t-\nu_t=\Delta^3h_Z(t)$),
where $\Delta$ is the difference operator. For any sequence $v=(v_1,v_2,\ldots,v_r)$
we define $\Delta(v)$ to be $(v_1,v_2-v_1,\ldots, v_r-v_{r-1})$. For a function 
$f$ defined on $\{0,1,2,\ldots\}$, we will define $\Delta(f)$ to be the function
$\Delta(f)(0)=f(0)$ and $\Delta(f)(t)=f(t)-f(t-1)$ for $t>0$.

\begin{defn}
Given $A = a_1p_1 + \cdots + a_rp_r \subset \pr N$, let $L$ be a hyperplane and let $F$ be the linear form defining $L$.
We do not assume $L$ passes through any point of $A$. The \defining{subscheme $B$ of $A$ residual to $L$}
is defined by the ideal $I_B=(I_A:F)$.
\end{defn}

Note that the ideal $(I_A:F)$ is saturated and the subscheme $B$ it defines is the fat point subscheme
$B = b_1 p_1 + \dots + b_r p_r$ where $b_i = a_i$ if $p_i \notin L$ and $b_i = \max(a_i-1,0)$ if $p_i \in L$.
We write $B = A:L$, as in
the following standard short exact sequence of graded 
$\field[\pr N]$-modules:
\[  0 \to I_{A:L}(-1) \overset{\cdot F}{\longrightarrow} I_A\to I_A/(I_A \cap (F)) \to 0 . \]
Sheafifying the above short exact sequence and tensoring by $\O_{\pr N}(t)$ yields the following short exact sequence
of $\O_{\pr N}$-modules, which will play the role of \eqref{fundses}:

\addtocounter{thm}{1} 
\begin{equation}\label{eqn:residual-ses}
0 \to \I_{A:L}(t-1) \to \I_A(t) \to \I_{L\cap A, L}(t) \to 0 
\end{equation}
where $\I_{L\cap A, L}$ denotes the ideal sheaf of
$L\cap A$ regarded as a subscheme of $L$ (or, depending on context, its extension by zero to $\pr N$).
When $N=2$, then $L\cong\pr1$ and $\I_{L\cap A, L}(t)\cong \O_{\pr1}(t-d)$, where
$d = \sum \{ a_i : p_i \in L \}$ is the degree of $L\cap A$.

We will use the well-known fact that $H^0(\I_A(t))$ can be identified with the homogeneous component 
$(I_A)_t$ of $I_A$ in degree $t$ (see, for example, \cite[Proposition 4.1.1]{refKrakow}).

\begin{remark}
Since \eqref{eqn:residual-ses} is the critical technical tool of this paper,
we provide here a brief justification.
It is easy to check that $I_{A:L}=I_B$ and $I_A$ sheafify to $\I_B$ and $\I_A$ respectively.
The inclusion $FI_B \subset I_A$ sheafifies to $\I_B\otimes \O_{\pr N}(-L)\subset \I_A$,
which is conventionally written as $\I_B(-1)\subset \I_A$.
Also, $I_A/FI_B = I_A/(I_A \cap (F)) = (I_A+(F)) / (F)$ is an ideal in
the homogeneous coordinate ring $R/(F)=K[L]\cong K[\pr {N-1}]$ of $\pr {N-1}$.
In particular $\I_A/\I_B(-1)$ is the sheafification of $I_A/FI_B = (I_A+(F))/(F)$,
so it is a subsheaf of the sheafification $\O_L$
of $K[L]$ (regarded as a sheaf on $\pr N$ by extension by 0).
Thus we have an exact sequence $0\to \I_B(-1)\to \I_A\to \I_{L\cap A, L}$.
By checking stalks, we see that in fact $\I_A\to \I_{L\cap A, L}$ is 
surjective. Thus $0\to \I_B(-1)\to \I_A\to \I_{L\cap A, L}\to 0$ is exact 
and $\I_{L\cap A, L}$ is the sheafification of $I_A/(I_A \cap (F))$.  
\end{remark}

By iteration we obtain from \eqref{eqn:residual-ses} the linked short exact sequences that will play the role of
\eqref{fundexactseqs}.

\begin{defn}
Let $A = a_1p_1 + a_2p_2 + \cdots + a_rp_r$ be a fat point scheme in $\pr N$.  Fix a sequence $L_1, \dots, L_{n+1}$ 
of hyperplanes in $\pr N$, not necessarily distinct.
\begin{enumerate}
\item[(a)] We define fat point schemes $A_0, \dots, A_{n+1}$ by $A_0 = A$ and 
$A_j = A_{j-1}:L_j$ for $1\le j\le n+1$.  
As a matter of notation, we write $A:L_1, \ldots,  L_j = (((A:L_1):L_2): \cdots : L_j)$, so we have $A_j = A:L_1, \ldots, L_j$.
\item[(b)] We will say that the sequence $L_1,\dots,L_{n+1}$ \defining{totally reduces} 
$A$ if $A_{n+1} = \varnothing$ is the empty scheme, or equivalently, if for each point $p_i$,
the hyperplane $L_j$ passes through $p_i$ for at least $a_i$ of the indices $j$.
\item[(c)] We associate to $A$ and the sequence $L_1,\dots,L_{n+1}$ an integer vector 
$\dvec=\dvec(A;L_1,\dots,L_{n+1})=(d_1,\dots,d_{n+1})$,
where $d_j$ is the degree $\deg(L_j\cap A_{j-1})$ of the scheme theoretic intersection of $L_j$ with $A_{j-1}$.
We refer to $\dvec$ as the \defining{reduction vector} for $A$ induced by
the sequence $L_1, \dots, L_{n+1}$. We will say that $\dvec$ is a \defining{full} reduction vector for $A$
if $L_1, \ldots, L_{n+1}$ totally reduces $A$.
\end{enumerate}
\end{defn}

\begin{remark}\label{degrem}
If $N=2$ and $L_1, \ldots, L_{n+1}$ totally reduce $A$ with reduction vector $\dvec = (d_1,\ldots,d_{n+1})$,
we note that $d_j$ is the sum of the multiplicities $a_{j-1,i}$ of $A_{j-1}=a_{j-1,1}p_1+\dots+a_{j-1,r}p_r$ 
at those points $p_i$ of $A_{j-1}$ which lie on $L_j$. More generally it is not hard to see for any $N$ that
$\deg(A_i) = d_{i+1} + \dotsb + d_{n+1}$ for $0 \leq i \leq n$ and
$\deg(A) = d_{1} + \cdots + d_{n+1}$. 
\end{remark}

\section{Bounds for the Hilbert function}\label{section2}

In this section we obtain various versions of our bounds on Hilbert functions
and give several examples.

\subsection{Bounds for the Hilbert function of the ideal}\label{sectionBnds}

In the context of hyperplanes $L_1,\ldots,L_{n+1}$ which totally reduce $A = a_1p_1 + a_2p_2 + \cdots + a_rp_r\subset \pr N$, let $\overline{A_i}$ be the scheme theoretic intersection $L_{i+1}\cap A_i$. Then
\eqref{fundexactseqs} becomes:
\addtocounter{thm}{1}
\begin{equation}\label{P2ses}
\begin{array}{@{0 \hspace{\tabcolsep} \to \hspace{\tabcolsep} } c @{\hspace{\tabcolsep} \to 
\hspace{\tabcolsep}} c @{\hspace{\tabcolsep} \to \hspace{\tabcolsep}} c @{\hspace{\tabcolsep} \to \hspace{\tabcolsep} 0}}
 \I_{A_{1}}(t-1) & \I_{A_0}(t) & \I_{\overline{A_0}, L_1}(t) \\
 \I_{A_{2}}(t-2) & \I_{A_1}(t-1) & \I_{\overline{A_1},L_2}(t-1) \\
 \multicolumn{3}{c}{\vdots} \\
 \I_{A_{i+1}}(t-(i+1)) & \I_{A_{i}}(t-i) & \I_{\overline{A_i},L_{i+1}}(t-i) \\
 \multicolumn{3}{c}{\vdots} \\
 \I_{A_n}(t-n) & \I_{A_{n-1}}(t-(n-1)) & \I_{\overline{A_{n-1}},L_n}(t-(n-1)) \\
 \O_{\pr N}(t-(n+1)) & \I_{A_n}(t-n) & \I_{\overline{A_n}, L_{n+1}}(t-n)
\end{array}
\end{equation}

From Theorem \ref{genthm} applied to \eqref{P2ses}, we almost immediately obtain the following result:

\begin{cor}\label{pNbnds}
Let $A_0=\sum_i m_i p_i$ be a fat point subscheme of $\pr N$ and let
$L_1,\ldots,L_{n+1}$ be hyperplanes which totally reduce $A_0$.
Let $h_{n}=h^0(L_{n+1},\I_{\overline{A_{n}},L_{n+1}}(t-n))+\binom{t-(n+1)+N}{N}$, for $0\leq i< n$ let
\[
  h_i=h^0(L_{i+1},\I_{\overline{A_i},L_{i+1}}(t-i))+h^0(\O_{\pr N}(t-(i+1)))-\deg(A_{i+1}),
\]
let $h_{n+1}'=h^0(\O_{\pr N}(t-(n+1)))$ and for $0\leq i\leq n$ let
\[
  h_i'=h^0(\O_{\pr N}(t-i))-\deg(A_i).
\]
Then 
\[
  \max_{0\leq i\leq n+1}(h_i') \leq \max_{0\leq i\leq n}(h_i)
    \leq h^0(\pr N,\I_{A_0}(t))
    \leq \binom{t-n+N-1}{N} + \sum_{i=0}^n h^0(L_{i+1},\I_{\overline{A_i},L_{i+1}}(t-i)).
\]
\end{cor}

\begin{proof}
The upper bound is immediate from Theorem \ref{genthm}
using the fact that $h^0(\pr N, \O_{\pr N}(t-(n+1)))=\binom{t-(n+1)+N}{N}$.
The lower bound $\max_{0\leq i\leq n}(h_i)\leq h^0(\pr N,\I_{A_0}(t))$
is also precisely the bound given in Theorem \ref{genthm}. To see this,
use the fact that $h^1(\pr N,\O_{\pr N}(t-(n+1)))=0$ for all $t$, so 
\[
  \max\big(h^0(\O_{\pr N}(t-(n+1))), h^0(L_{n+1},\I_{\overline{A_{n}},L_{n+1}}(t-n))+h^0(\O_{\pr N}(t-(n+1)))-h^1(\O_{\pr N}(t-(n+1)))\big)
\]
is $h_{n}=h^0(L_{n+1},\I_{\overline{A_{n}},L_{n+1}}(t-n))+\binom{t-(n+1)+N}{N}$. The remaining terms 
over which the maximum is taken in 
\eqref{eqn: ideal lower bounds} are of the form
\addtocounter{thm}{1}
\begin{equation}\label{oneterm}
h^0(L_{i+1},\I_{\overline{A_i},L_{i+1}}(t-i))+h^0(\I_{A_{i+1}}(t-(i+1)))-h^1(\I_{A_{i+1}}(t-(i+1)))
\end{equation}
for $0\leq i< n$. Taking cohomology of the short exact sequence
\[  0 \to \I_{A_{i+1}}(t-(i+1)) \to \O_{\pr N}(t-(i+1)) \to \O_{A_{i+1}}(t-(i+1)) \to 0  \]
and using $h^1(\O_{\pr N}(t-(i+1)))=0$ gives 
\[  h^0(\I_{A_{i+1}}(t-(i+1))) - h^1(\I_{A_{i+1}}(t-(i+1))) = h^0(\O_{\pr N}(t-(i+1))) - h^0(\O_{A_{i+1}}(t-(i+1))),  \]
but $h^0(\O_{A_{i+1}}(t-(i+1))) = h^0(\O_{A_{i+1}}) = \deg(A_{i+1})$ since $\O_{A_{i+1}}(t-(i+1))\cong \O_{A_{i+1}}$ 
(because $\O_{A_{i+1}}$ is a skyscraper sheaf). Thus \eqref{oneterm} becomes
\[  h_i=h^0(L_{i+1},\I_{\overline{A_i},L_{i+1}}(t-i))+h^0(\O_{\pr N}(t-(i+1)))-\deg(A_{i+1})  \]
which is exactly what occurs in the lower bound $\max_{0\leq i\leq n}(h_i)$.
Moreover, from the long exact sequence
of cohomology coming from the top sequence in \eqref{P2ses} we see that
\eqref{oneterm} with $i=0$ (which we just saw is $h_0$) is at least as big as
$h^0(\I_{A_0}(t))-h^1(\I_{A_0}(t))$, which by the same argument as above is $h_0'$.
Thus $h_0\geq h_0'$ and clearly $h_i\geq h_{i+1}'$ for $0\leq i\leq n$,
so we now see that $\max_{0\leq i\leq n+1}(h_i')\leq\max_{0\leq i\leq n}(h_i)$.
\end{proof}

\begin{rem}\label{refrem1}
When $N>2$ we will not in general know the terms $h^0(\I_{\overline{A_i},L_{i+1}}(t-i))$, and thus we 
will not in general know either the upper bound in Corollary \ref{pNbnds} or the lower bound
$\max_{0\leq i\leq n}(h_i)$. Nevertheless, by induction on dimension
we can obtain upper and lower bounds on $h^0(\I_{\overline{A_i},L_{i+1}}(t-i))$ and 
thus also on $h^0(\I_{{A_0}}(t))$.
On the other hand, we always will know each $h_i'$ and thus we know the lower bound 
$\max_{0\leq i\leq n+1}(h_i')$. Easy examples (such as three collinear points 
of multiplicity 1 in $\pr3$) show in general that equality in
$\max_{0\leq i\leq n+1}(h_i')\leq \max_{0\leq i\leq n}(h_i)$
fails. When $N=2$, however, equality holds and we have 
$h^0(\I_{\overline{A_i},L_{i+1}}(t-i))=\binom{t-i-d_{i+1}+1}{1}$, so our bounds for $N=2$, 
which we state explicitly in the next corollary,
can be computed exactly. (To see that equality holds when $N=2$,
note that $h_i=h_{i+1}'$ if $t-i-d_{i+1}+1< 1$, while
if $t-i-d_{i+1}+1\geq 1$, then $h^1(\I_{\overline{A_i},L_{i+1}}(t-i))=0$, so $h_i= h_i'$.)
\end{rem}

\begin{cor}\label{idealbounds}
Let $A_0=\sum_i m_i p_i$ be a fat point subscheme of $\pr2$, let $L_1,\ldots,L_{n+1}$ 
be a totally reducing sequence of lines with reduction vector
$\dvec = (d_1,\ldots,d_{n+1})$, let $h_{n+1}'= \binom{t-(n+1)+2}{2}$
and for $0\leq i\leq n$ let 
\[  h_i'=\binom{t-i+2}{2}-\sum_{i+1\leq j\leq n+1}d_j.  \]
Then we have
\[
  \max(h_0',\ldots,h_{n+1}') \leq h^0(\pr2,\I_{A_0}(t)) \leq \binom{t-n+1}{2} + \sum_{i=0}^n\binom{t-i-d_{i+1}+1}{1}.
\]
\end{cor}

\begin{proof}
This follows immediately from Corollary \ref{pNbnds}
using the facts that $\deg(A_i)=d_{i+1}+\cdots+d_{n+1}$ and
$h^0(L_{i+1},\I_{\overline{A_{i}},L_{i+1}}(t-i))=h^0(\pr1,\O_{\pr1}(t-i-\deg(\overline{A_i})))
=\binom{t-i-d_{i+1}+1}{1}$. 
\end{proof}

\begin{rem}\label{refrem} In Corollary \ref{idealbounds}, the referee notes
that the lower bounds $h_i'\leq h^0(\pr2,\I_{A_0}(t))$ can be seen very simply
using the obvious facts that $h_i'=\binom{t-i+2}{2}-\deg(A_i)\leq h^0(\pr2,\I_{A_i}(t-i))$
and $h^0(\pr2,\I_{A_i}(t-i))\leq h^0(\pr2,\I_{A_0}(t))$.
Moreover, in the case that the lines $L_i$ are distinct,
the referee also shared with us a nice way of obtaining the upper bound. Let $Z'$ be the fat point scheme
resulting from adding $\binom{t-i-d_{i+1}+1}{1}$ points (in general position on line $L_i$) to $Z$,
let $Z_{n+1}$ and $Z'_{n+1}$ be the residuals of $Z$ and $Z'$ respectively with respect to
the sequence of lines $L_i$, and note that $Z'_{n+1}=Z_{n+1}$.
The added points impose at most $\sum_i \binom{t-i-d_{i+1}+1}{1}$ additional conditions to forms of
degree $t$ vanishing on $Z$, hence $h_{I_Z}(t)\leq h_{I_{Z'}}(t)+\sum_i \binom{t-i-d_{i+1}+1}{1}$.
But the number of points added to each line is such that each line is forced by B\'ezout's Theorem to be
a fixed component of $(I_{Z'})_t$, so $h_{I_{Z'}}(t)=h_{I_{Z'_{n+1}}}(t-(n+1))$, and since the
lines are totally reducing we have that $Z'_{n+1}=Z_{n+1}$ is empty, so $h_{I_{Z'_{n+1}}}(t-(n+1))=\binom{t-(n+1)+2}{2}$.
\end{rem}

\subsection{Equality of the bounds}\label{subsectionEquality}

We now define a property that implies the bounds in Corollary \ref{idealbounds} coincide.

\begin{defn}\label{defn: GMS}
A vector $\vvec = (v_1, \ldots, v_r)$ will be said to be \defining{non-negative},
respectively \defining{positive}, \defining{non-increasing} or \defining{strictly decreasing} if its entries are.
If $\vvec=(v_1,\dots,v_r)$ is a non-negative, non-increasing integer vector, 
we will say $\vvec$ is \defining{\GMS} if $v_i - v_j \geq j-i-1$ for each $1 \leq i < j \leq r$.
\end{defn}

\begin{thm}\label{GMSthm}
Let $\dvec = (d_1,\ldots,d_{n+1})$ be the reduction vector
of a fat point subscheme of $\pr2$ with respect to some 
totally reducing sequence of lines $L_1,\ldots,L_{n+1}$.
If $\dvec$ is \GMS, then the upper and lower bounds given
in Corollary \ref{idealbounds} coincide.
\end{thm}

\begin{proof} By Remark \ref{equalityrem} it is enough to show that 
$h^0(\O_{L_i}(t-(i-1)-d_i))>0$ implies $h^1(\O_{L_j}(t-(j-1)-d_j))=0$
for all $j>i$ for all $t$. I.e., it suffices to show that
$t-(i-1)-d_i\ge0$ implies $t-(j-1)-d_j\ge-1$,
but \GMS\ means $d_i-d_j-j+i+1\ge 0$, and adding
$t-(i-1)-d_i\ge0$ to $d_i-d_j-j+i+1\ge 0$ gives $t-(j-1)-d_j\ge-1$.
\end{proof}

The converse of Theorem \ref{GMSthm} usually also holds; see Remark \ref{thefunctions2}.

\begin{remark}\label{GMSversionofbezout}
Consider $\vvec = (v_1, \ldots, v_r)$.
The \GMS\ property is weaker than being strictly monotone (each $v_i - v_j \geq j-i$),
but stronger than being non-increasing (each $v_i-v_j \geq 0$).
The adjective ``\GMS'' is meant to suggest ``generalized monotone sequence''.
In \cite{refGMS}, a structure analogous to our reduction vector
is called a pseudo type vector. When a pseudo type vector
satisfies a certain property (specifically Property 3.2 in the statement of Theorem 3.7 of \cite{refGMS}), 
\cite{refGMS} shows how to build a double point scheme
whose Hilbert function is determined by the pseudo type vector.
They build it starting from the empty scheme in a series of steps
(called basic double links) related to the entries
of the pseudo type vector. We on the other hand start with an arbitrary fat point scheme
and in a series of steps corresponding to the entries
of a reduction vector reduce the scheme to the empty scheme;
when the reduction vector is \GMS, it uniquely determines the Hilbert function of the scheme.
Because we tear down where \cite{refGMS} builds up,
given a double point scheme for which we have a reduction vector $\vvec = (v_1, \ldots, v_r)$,
the corresponding pseudo type vector would be $(v_r, \ldots, v_1)$
(which is the ordering we used in previous versions of our paper also).
Property 3.2 of \cite{refGMS} applied to $(v_r, \ldots, v_1)$ is 
precisely that between any two zero entries of 
$\diff (v_r, \ldots, v_1)$ there is an entry strictly bigger than 1.
We will show in Proposition \ref{GMSprop} that this
is equivalent to our condition \GMS, which thus, if you like, can be
taken to stand for the authors ``Geramita, Migliore, Sabourin'' of \cite{refGMS}.

It is also true that when we have a \GMS\ reduction vector for a fat point scheme
$Z$, the residuation steps taken in reverse are actually basic double links.
Thus our results show how to construct fat point schemes other than just double points as
a sequence of basic double links. We direct readers interested in further details 
to section 7.2 of our expository posting \cite{refCHT}. 
\end{remark}

\begin{prop}\label{GMSprop}
Let $\vvec = (v_1, \ldots, v_r)$ be a non-negative, non-increasing integer vector.
Then the following are equivalent:
\begin{itemize}
\item[(a)] $\vvec$ is \GMS;
\item[(b)] between any two zero entries of 
$\diff (v_r, \ldots, v_1)$ there is an entry strictly bigger than 1;
\item[(c)] $\vvec$ does not contain a subsequence of consecutive entries of the form $(a,a,a)$,
or of the form $(a_i,\ldots,a_{i+j+1})$ for $j>1$ where $a_i=a_{i+1}$, $a_{i+j}=a_{i+j+1}$,
and $a_{i+1},\ldots,a_{i+j}$ are consecutive integers.
\end{itemize}
\end{prop}

\begin{proof} Assume (a) holds. Then (b) holds by the pigeon hole principle. To show (b) implies (c),
it is enough to check the contrapositive, which is trivial. To show that (c) implies (a),
apply the pigeon hole principle to the contrapositive.
\end{proof}

\subsection{Bounds on the Hilbert function of the scheme}\label{bndsonschemesubsect}

We can reformulate the bounds of Corollary \ref{idealbounds} in terms of the Hilbert function of $A_0$, 
rather than in terms of
the ideal of $A_0$. The lower bound of the theorem becomes the upper bound 
\addtocounter{thm}{1}
\begin{equation}\label{tag1}
 h_{A_0}(t) \leq \min_{0\leq i \leq n+1} \Big( h^0(\pr2,\O_{\pr2}(t)) -  h_i'  \Big) 
 = \min_{0\leq i \leq n+1} \Bigg( \binom{t+2}{2} - \binom{t-i+2}{2} + \sum_{i+1\leq j\leq n+1} d_j \Bigg).
\end{equation}

From the upper bound
$h^0(\pr2,\I_{A_0}(t)) \leq \binom{t-n+1}{2} + \sum_{i=0}^n\binom{t-i-d_{i+1}+1}{1}$ we obtain
\begin{equation*}
\begin{split}
h_{A_0}(t) & \geq  \binom{t+2}{2}- \binom{t-n+1}{2} - \sum_{i=0}^n\binom{t-i-d_{i+1}+1}{1}\\
& = \sum_{i=0}^n \Bigg(\binom{t-i+1}{1} - \binom{t-i-d_{i+1}+1}{1} \Bigg).
\end{split}
\end{equation*}

We can simplify each term in this last sum with the following identity:
if $a \geq b$ then,
by examining the three cases $0 \leq b \leq a$, $b \leq 0 \leq a$, and $b \leq a \leq 0$,
we obtain
\[
 \binom{a}{1} - \binom{b}{1} =\left.
   \begin{cases}
     a-b , & 0 \leq b \leq a \\
     a    , & b \leq 0 \leq a \\
     0    , & b \leq a \leq 0
   \end{cases}  
   \right\rbrace = \binom{\min(a,a-b)}{1}.
\]
We can apply this to each term in our lower bound.
This yields
\addtocounter{thm}{1}
\begin{equation}\label{tag3}
h_{A_0}(t) \geq \sum_{i=0}^n \binom{\min(t-i+1 , d_{i+1})}{1}.
\end{equation}

\begin{defn}\label{thefunctions1} 
Given an integer vector $\vvec = (v_1, \ldots, v_{n+1})$, we define functions
\[  f_\vvec(t)=\sum_{i=0}^n \binom{\min(t-i+1 , v_{i+1})}{1}  \]
and 
\[  F_\vvec(t)=\min_{0 \leq i \leq n+1} \Bigg( \binom{t+2}{2} - \binom{t-i+2}{2} + \sum_{i+1\leq j\leq n+1}v_j\Bigg).  \]
\end{defn}

\subsection{Optimality issues}
An obvious question is how in Corollary \ref{idealbounds} should one pick the lines $L_i$ in order to get the tightest
bounds for a given fat point subscheme $A\subset\pr2$. This is not clear. 
Experience suggests that one should not deviate too much from a greedy approach in which one tries to pick
$L_{i+1}$ so that $d_{i+1}=\deg(L_{i+1}\cap A_i)$ is as big as possible.
For example, one can always avoid choices such that $\dvec$ fails to be non-increasing, which one must do if
one hopes to have a \GMS\ reduction vector and thereby determine $h_A$ exactly by an application of Theorem \ref{GMSthm}. 

The next example shows, for our lower bound $f_\dvec$ at least,
that making greedy choices to obtain $\dvec$ does not always give the best bound. 
This example also shows that by choosing
different sequences, our bounds can sometimes determine the Hilbert function exactly,
even when neither sequence is \GMS.

\begin{example}\label{greedydetermine}
Consider the scheme $A$ consisting of three points of multiplicity 2
and twelve points of multiplicity 1
arranged as shown in Figure~\ref{fig:nongreedy}.

\addtocounter{thm}{1}
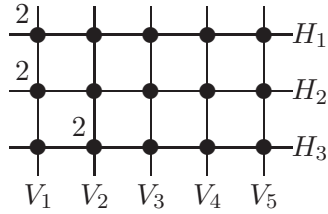
\begin{figure}[h]
\setlength{\unitlength}{0.75cm}
\begin{picture}(5,4)(0,-1)
\multiput(0,0)(1,0){5}{\circle*{.3}}
\multiput(0,1)(1,0){5}{\circle*{.3}}
\multiput(0,2)(1,0){5}{\circle*{.3}}
\put(-0.4,1.2){$2$}
\put(-0.4,2.2){$2$}
\put(0.6,0.2){$2$}
\multiput(-.5,0)(0,1){3}{\line(1,0){5}}
\multiput(0,-.5)(1,0){5}{\line(0,1){3}}
\put(4.5,1.85){$H_1$}
\put(4.5,.85){$H_2$}
\put(4.5,-.15){$H_3$}
\put(-.25,-1){$V_1$}
\put(.75,-1){$V_2$}
\put(1.75,-1){$V_3$}
\put(2.75,-1){$V_4$}
\put(3.75,-1){$V_5$}
\end{picture}
\caption{A configuration where greedy choices are not optimal.}\label{fig:nongreedy}
\end{figure}

The sequence of lines $L_1, \ldots, L_5$, where
$L_1=H_1$, 
$L_2=H_2$, 
$L_3=H_3$, 
$L_4=V_1$, and
$L_5=V_2$ 
is totally reducing for 
$A$
and gives 
$\dvec=(6,6,6,2,1)$ 
as the reduction vector. The sequence of lines $L_1', \ldots, L_7'$, where
$L'_1=V_1$, 
$L'_2=V_2$, 
$L'_3=V_3$, 
$L'_4=V_4$, 
$L'_5=V_5$,
$L'_6=V_1$, and
$L'_7=V_2$ 
is also totally reducing for 
$A$
and gives $\dvec'=(5,4,3,3,3,2,1)$ as the reduction vector.

Note that $\dvec$ is obtained by making only greedy choices, whereas $\dvec'$ is not.
Neither $\dvec$ nor $\dvec'$ is \GMS, so $f_\dvec \neq F_\dvec$ and $f_{\dvec'} \neq F_{\dvec'}$.
Hence neither determines $h_A$.
Nevertheless we compute
\[
  \begin{split}
  f_\dvec &= (1, 3, 6, 10, 15, 18, 20, 21,21,\dots) \\
  F_\dvec &= (1, 3, 6, 10, 15, 18, 21,21,\dots) \\
  f_{\dvec'} &= (1, 3, 6, 10, 15, 18, 21,21,\dots) \\
  F_{\dvec'} &= (1, 3, 6, 10, 15, 21, 21,\dots)
  \end{split}
\]
so that $h_A \leq F_\dvec = f_{\dvec'} \leq h_A$, which does determine $h_A$ uniquely.
Observe that $f_\dvec(6) = 20 < 21 = f_{\dvec'}(6)$, so $\dvec'$ gives a strictly better lower bound
than the reduction vector $\dvec$ obtained by making greedy choices.
\end{example}

\begin{example}\label{ZachExample}
In this example, we have a subscheme $A=2p_1+p_2+p_3+2p_4+p_5+p_6$
as shown in Figure \ref{example:  intro ex2}
and we are given information about which points are collinear.
In particular, for each line $\ell_i$ we indicate by a 1 each point $p_j$
which lies on $\ell_i$, and by 0 if $p_j$ does not lie on $\ell_i$. We show
the bounds we obtain on $h_Z$ in this case are best possible with the given information.  

\addtocounter{thm}{1}
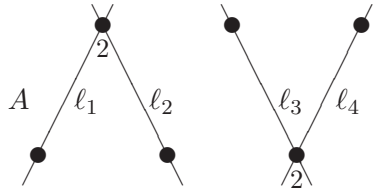
\begin{figure}[h]
\begin{tabular*}{0.9\textwidth}{c@{\extracolsep{\fill}}c}
\begin{minipage}[c]{0.5\textwidth}
 \begin{tabular}{c|cccccc}
 Points   & $p_1$ & $p_2$ & $p_3$ & $p_4$ & $p_5$ & $p_6$ \\
 \hline 
 Multiplicities & 2 & 1 & 1 & 2 & 1 & 1 \\
 \hline
 $\ell_1$ & 1 & 1 & 0 & 0 & 0 & 0 \\
 $\ell_2$ & 1 & 0 & 1 & 0 & 0 & 0 \\
 $\ell_3$ & 0 & 0 & 0 & 1 & 1 & 0 \\
 $\ell_4$ & 0 & 0 & 0 & 1 & 0 & 1
 \end{tabular}
\end{minipage}
&
\begin{minipage}[c]{0.35\textwidth}
\setlength{\unitlength}{0.4cm}
\begin{picture}(4,4)(0,0)
\put(2.5,5){\line(1,-2){3}}
\put(.2,-1){\line(1,2){3}}
\put(2.85,4.3){\circle*{0.5625}} 
\put(2.63,3.25){\small 2}
\put(.7,0){\circle*{0.5625}}
\put(-0.3,1.5){$A$}
\put(1.9,1.5){$\ell_1$}
\put(4.4,1.5){$\ell_2$}
\put(8.7,1.5){$\ell_3$}
\put(10.6,1.5){$\ell_4$}
\put(5,0){\circle*{0.5625}}
\put(6.8,5){\line(1,-2){3}}
\put(8.8,-1){\line(1,2){3}}
\put(7.15,4.3){\circle*{0.5625}}
\put(9.3,0){\circle*{0.5625}}  
\put(9.08,-1.1){\small 2}
\put(11.45,4.3){\circle*{0.5625}}
\end{picture}
\end{minipage}
\end{tabular*} 
\caption{Collinearity data for a given fat point scheme.}\label{example:  intro ex2}
\end{figure}
One easily checks that $L_1 = \ell_1, L_2 = \ell_3, L_3 = \ell_2, L_4 = \ell_4$ is a sequence of lines 
yielding the reduction vector $\dvec=(3,3,2,2)$.  
We then find
\[
\begin{split}
 f_{(3,3,2,2)} &= (1,3,6,9,10,10, \ldots), \\
 F_{(3,3,2,2)} &= (1,3,6,10,10, \ldots).
\end{split}
\]
\addtocounter{thm}{1}
\begin{figure}[h]
\vskip\baselineskip
\begin{tabular*}{0.9\textwidth}{c@{\extracolsep{\fill}}c}
\begin{minipage}[c]{0.45\textwidth}
\setlength{\unitlength}{0.4cm}
\begin{picture}(4,4)(-4,0)
\put(2.5,5){\line(1,-2){3}}
\put(.2,-1){\line(1,2){3}}
\put(2.85,4.3){\circle*{0.5625}} 
\put(2.63,3.25){\small 2}
\put(.7,0){\circle*{0.5625}}
\put(-0.3,1.5){$A'$}
\put(1.9,1.5){$\ell_1$}
\put(4.4,1.5){$\ell_2$}
\put(8.7,1.5){$\ell_3$}
\put(10.6,1.5){$\ell_4$}
\put(5,0){\circle*{0.5625}}
\put(6.8,5){\line(1,-2){3}}
\put(8.8,-1){\line(1,2){3}}
\put(7.15,4.3){\circle*{0.5625}}
\put(9.3,0){\circle*{0.5625}}  
\put(9.08,-1.1){\small 2}
\put(11.45,4.3){\circle*{0.5625}}
\multiput(0,4.3)(.3,0){44}{\circle*{0.1}}
\multiput(0,0)(.3,0){44}{\circle*{0.1}}
\multiput(1.8825,4.915)(.3225,-0.215){27}{\circle*{0.1}}
\end{picture}
\end{minipage}
&
\begin{minipage}[c]{0.45\textwidth}
\setlength{\unitlength}{0.4cm}
\begin{picture}(4,4)(-3,-.1)
\put(2.5,5){\line(1,-2){3}}
\put(.2,-1){\line(1,2){3}}
\put(2.85,4.3){\circle*{0.5625}}
\put(2.63,3.25){\small 2}
\put(.7,0){\circle*{0.5625}}
\put(-0.3,1.5){$A''$}
\put(1.9,1.5){$\ell_1$}
\put(4.4,1.5){$\ell_2$}
\put(8.7,1.5){$\ell_3$}
\put(10.6,1.5){$\ell_4$}
\put(4.65,.65){\circle*{0.5625}}
\put(6.8,5){\line(1,-2){3}}
\put(8.8,-1){\line(1,2){3}}
\put(7.15,4.3){\circle*{0.5625}}
\put(9.3,0){\circle*{0.5625}}
\put(9.08,-1.1){\small 2}
\put(11.45,4.3){\circle*{0.5625}}
\multiput(0,4.3)(.3,0){44}{\circle*{0.1}}
\multiput(0,0)(.3,0){44}{\circle*{0.1}}
\multiput(1.8825,4.915)(.3225,-0.215){27}{\circle*{0.1}}
\end{picture}
\end{minipage}
\end{tabular*}
\vskip\baselineskip
\caption{Two fat point schemes compatible with the data in Figure \ref{example:  intro ex2}.}\label{example:  intro ex2'}
\end{figure}
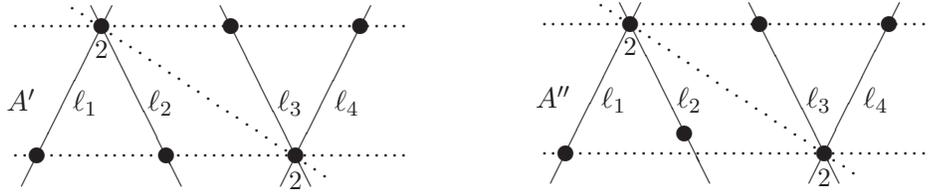
We now show that the bounds $f_{(3,3,2,2)}$ and $F_{(3,3,2,2)}$ are sharp.
Consider the schemes $A'$ and $A''$ of Figure \ref{example:  intro ex2'}; these are compatible with the collinearity data given in
Figure \ref{example:  intro ex2}, but now we are given additional collinearity data, represented by the dotted lines.
In particular, there exists a reduction vector $\vvec'=(4,4,2)$ for $A'$ obtained using the dotted lines. In this case
we have $f_{\vvec'}=F_{\vvec'}=(1,3,6,9,10,10,\ldots)$ so $h_{A'}(3)=9$.
For the fat point scheme $A''$ in Figure \ref{example:  intro ex2'} one can choose lines to obtain a reduction vector
$\vvec''=(4,3,2,1)$, giving $f_{\vvec''}=F_{\vvec''}=(1,3,6,10,10,\ldots)$, hence $h_{A''}(3)=10$.
Since both $f_{(3,3,2,2)}$ and $F_{(3,3,2,2)}$ are Hilbert functions of fat point schemes
having reduction vector $(3,3,2,2)$, the bounds $f_{(3,3,2,2)}$ and $F_{(3,3,2,2)}$ are optimal.

\end{example}

\subsection{Further remarks}

\begin{remark}\label{thefunctions2}
When $Z$ is a fat point subscheme of $\pr2$ for which $\vvec$ occurs as 
a full reduction vector, then  
we have $f_\vvec(t)\leq h_Z(t)\leq F_\vvec(t)$ for all $t\ge 0$ (by 
Corollary \ref{idealbounds}, Equations \eqref{tag1} and \eqref{tag3}, and Definition \ref{thefunctions1}),
and if $\vvec$ is \GMS\ then by Theorem \ref{GMSthm} we have $f_\vvec=F_\vvec$. 
For any non-negative integer vector $\vvec$, this means that $f_\vvec(t)\leq F_\vvec(t)$ for all $t\geq0$,
with equality if $\vvec$ is \GMS\
(because any non-negative integer vector $\vvec=(v_1,\ldots,v_r)$ is a full reduction vector
for some $Z$; just take $Z$ to be the reduced scheme consisting of $v_i$ general points
on each line $L_i$ of a general set of lines $L_1,\ldots,L_r$). 
Conversely, if $\vvec$ is positive and if $f_\vvec=F_\vvec$, then it is possible to show
that $\vvec$ is \GMS\ (see \cite[Theorem 8.1.1]{refCHT}). In most cases
reduction vectors are positive (this is because a zero entry in a reduction vector means
that one has taken a residual with respect to a line
that doesn't meet the scheme; this is a null operation, so
one would normally avoid such lines). Thus equality of our bounds
$f_\vvec$ and $F_\vvec$ is essentially equivalent to $\vvec$ being \GMS.
\end{remark}

\begin{remark}\label{thefunctions3} 
For those interested in writing code, recursive definitions for the functions $f_\dvec$ and $F_\dvec$
may be of interest. The recursions can be justified by checking that the explicit versions of the functions
shown above satisfy the recursion. Here we merely exhibit the result; we refer the interested reader to \cite{refCHT} 
for details.

Let $\vvec=(v_1,\ldots,v_r)$ be a  non-negative integer vector. Then $f_\vvec$ can be defined by
setting $f_\vvec(t)=0$ if $t < 0$ or $n=0$ (i.e., if $\vvec$ is the empty vector), and
otherwise we set $f_\vvec(t) =  f_{\vvec'}(t-1) + \min( t+1 , v_1 )$,
where $\vvec' = (v_2,\dots,v_r)$.
Similarly, if $t < 0$ or $\vvec$ is the empty vector, then $F_\vvec(t)=0$.
Otherwise we set
$F_\vvec(t) = \min \Big( t+1+F_{\vvec'}(t-1), \sum_{i=1}^r v_i \Big)$.
\end{remark}

One additional way to define $f_\vvec$ may be of interest, since it shows that
$f_\vvec$ is what \cite{refGMS} calls the {\it standard Hilbert function\/}
associated to $\vvec$. 

\begin{defn}\label{def:standardconfig}
Let $\vvec=(v_1,\dots,v_n)$ be a non-negative integer vector.
Define the \defining{summation operator} $\sum$ by
\[ \sum \vvec = (v_1, v_1+v_2, \dots, v_1 + \dots + v_n). \]
Also define the \defining{standard configuration} $S_\vvec$ determined by $\vvec$
to be the set $S_\vvec\subset {\bf Z}\times{\bf Z}$ 
of all integer lattice points $(i,j-1)$ with $i\ge0$ and $1\leq j\leq n$ such that $i<v_j$.
Thus $S_\vvec$ consists of the $v_j$ leftmost first quadrant lattice points
of ${\bf Z}\times{\bf Z}$ on each 
horizontal line with second coordinate $j-1$ for $1\leq j\leq n$
(we take $j-1$ to be the second coordinate so that the bottom row of
$S_\vvec$ is on the horizontal axis).  We also 
define the \defining{diagonal count operator} $\diag(\vvec)$ of $\vvec$
by $\diag(\vvec)=(c_0,c_1,\dots, c_t, \dots)$, where
$c_t$ is the number of points in $S_\vvec$ lying on the diagonal line
with the equation $x+y=t$.
Equivalently,
\[ c_t = \card\{\, j : 0 \leq j < n \text{ and } 0 \leq i < v_{j+1} \text{ where } i+j=t \, \} . \]
\end{defn}

\begin{remark}\label{thefunctions4} 
Given any non-negative integer vector $\vvec=(v_1,\dots,v_n)$, we have $f_\vvec = \sum \diag(\vvec)$.
Indeed, $\binom{\min(t-i+1 , v_i)}{1}$ is the number of dots in row $i$ of $S_\vvec$
(i.e., with vertical coordinate $y=i-1$) on or to the left of the line with slope $-1$
intersecting the horizontal coordinate axis at coordinate $x=t$.
Thus from Definition \ref{thefunctions1} we see that
$f_\vvec(t)$ is the number of points of the standard configuration $S_\vvec$ in the region $\{x+y \leq t\}$;
i.e., $f_\vvec=\sum\diag(\vvec)$, as claimed.
Note that this is precisely the standard Hilbert function
associated to $\vvec$ defined in \cite{refGMS}. See \cite{refCHT} for further discussion.
\end{remark}

\begin{example}
Let $\vvec=(8,6,5,2)$.
Then $S_\vvec$ is the set of lattice points in ${\bf Z}\times{\bf Z}$ shown in Figure \ref{fig:diagonals}.
We can regard $\vvec$ as giving the row counts for $S_\vvec$ and
$\diag(\vvec) = (1, 2, 3, 4, 4, 3, 3, 1, 0, 0, \dots)$
as giving
the diagonal counts (along diagonals with slope $-1$).

\addtocounter{thm}{1}
\setlength{\unitlength}{0.75cm}
\begin{figure}[h]
\begin{picture}(6,4.5)
\multiput(0,0)(1,0){8}{\circle*{.3}}
\multiput(0,1)(1,0){6}{\circle*{.3}}
\multiput(0,2)(1,0){5}{\circle*{.3}}
\multiput(0,3)(1,0){2}{\circle*{.3}}
\put(0,0){\line(1,0){8}}
\put(0,0){\line(0,1){4}}
\put(-0.1,4.4){$y$}
\put(8.4,-0.1){$x$}
\end{picture}
\caption{$S_{(8,6,5,2)}$.}\label{fig:diagonals}
\end{figure}
\end{example}

\begin{example}
Note that $\diag(\vvec)$ does not determine $\vvec$.
Indeed, $\diag((2,2)) = \diag((3,1)) = (1, 2, 1, 0, 0, \dots)$.
\end{example}

\begin{remark}
One may ask whether $f_\vvec$ and $F_\vvec$ are {\it differentiable O-sequences}, since
given a 0-dimensional subscheme $Z\subset\pr2$, $h_Z$ is a differentiable O-sequence
(for details, see \cite{refGMR, refMa, refS, refMi}).
This is equivalent to saying that $h_Z=\sum\diag(\dvec)$ for an integer vector $\dvec$
with positive increasing entries (see \cite[Corollary 3.4]{refGMR} and \cite[Theorem 2.6]{refGHS}).
In fact, if $\dvec = (d_1,\dots,d_r)$ is \GMS\ then
$f_{\dvec}$ and $F_{\dvec}$ are differentiable O-sequences, since by the construction
mentioned in Remark \ref{thefunctions2} every \GMS\ sequence occurs
as a full reduction vector for some fat point scheme $Z$, in which case
$f_\dvec=h_Z=F_{\dvec}$. More generally, if $\dvec$ is merely non-negative,
then $F_\dvec=h_Z$ for some fat point subscheme $Z$ 
(but $Z$ need not have $\dvec$ as a reduction vector),
and hence $F_\dvec$ is a differentiable O-sequence.
And if $\dvec$ is non-negative and non-increasing,
then $f_\dvec=h_{Z'}$ for some fat point subscheme $Z'$ (again not necessarily with $\dvec$ as a reduction vector),
and hence $f_\dvec$ is then also a differentiable O-sequence. See \cite[Remark 7.5.2]{refCHT}
for details.
\end{remark}

\section{Graded Betti numbers}\label{bettisection}

We focus in this section on the underlying principles that can be used to
give lower and upper bounds for the graded Betti numbers of a fat point scheme $A \subset \pr2$ 
in terms of a reduction vector $\dvec$. 

The graded Betti numbers
determine the modules in a minimal free graded resolution for $I_A$
up to graded isomorphism. When $h_A$ is known, the Betti numbers of the syzygy module
can be determined from the Betti numbers for the module of generators.
Thus we focus on bounding the Betti numbers $\nu_{t+1}$ for the module of generators,
given a \GMS\ reduction vector $\dvec$, where $\nu_{t+1}(A)$ is the number of generators
of $I_A$ of degree $t+1$ in any minimal set of homogeneous generators.
Explicitly, $\nu_{t+1}(A)$ is the dimension of the cokernel of the map
$\mu_t(A): (I_A)_t\otimes R_1\to (I_A)_{t+1}$
given by multiplication of forms in the ideal of degree $t$ by linear forms.

\subsection{The underlying principles}

Assuming that we have a full \GMS\ reduction vector $\dvec$ for a fat point subscheme
$A\subset\pr2$, then we know $h_A$ in terms of $\dvec$. In particular, 
we know $\alpha(I_A)$; i.e., the least degree $t\ge0$ such that $(I_A)_t\ne 0$. Hence we know that
$\nu_t(A)=0$ for $t<\alpha(I_A)$ and that $\nu_{\alpha(I_A)}=\dim (I_A)_{\alpha(I_A)}$.
We also know the Castelnuovo-Mumford regularity $\reg(I_A)$; i.e., the least degree
$t$ such that $h_A(t-1) = \deg A$. By the well-known fact that one can always choose
homogeneous generators for $I_A$ 
of degree at most $\reg(I_A)$, we see that $\nu_t(A)=0$ for $t>\reg(I_A)$.

It is easy to give naive bounds on $\nu_{t+1}(A)$ that depend only on $\dvec$ for degrees $t$ in the range 
$\alpha(I_A)\leq t < \reg(I_A)$. For example, we have 
\addtocounter{thm}{1}
\begin{equation}\label{naivebnds}
\max(0,h_{I_A}(t+1)-3h_{I_A}(t))\leq \nu_{t+1}(A)\leq h_{I_A}(t+1)-(2+h_{I_A}(t)).
\end{equation}
The lower bound comes from the fact that the dimension of the cokernel of a linear map of 
vector spaces is at least 0 but also at least the difference in the dimensions 
of the target and source of the map. Thus the lower bound here is an equality
if and only if $\mu_t$ has maximal rank.
The upper bound here is obtained
by applying the following lemma with $V=(I_A)_t\subset R_t$ for $R=\field[\pr 2]$:

\begin{lem}\label{lem:genimfact}
Let $t\ge0$ and let $V\neq0$ be a vector subspace of $R_t$ for $R=\field[\pr N]=\field[x_0,\dots,x_N]$. 
Then $R_1V\subseteq R_{t+1}$ has dimension $\dim_\field R_1V\geq N + \dim_\field V$.
\end{lem}

\begin{proof}Let $W_i\subseteq R_1$ be the $\field$-span of $x_0,\dots,x_i$. Then
$x_iV\not\subset W_{i-1}V$, since the highest degree in the $x_i$ variable among monomial
terms of elements is the same for $V$ as for $W_{i-1}V$, but it is higher for $x_iV$ than it is for $V$.  
Thus we have proper containments $W_0V\subsetneq W_1V\subsetneq \dots \subsetneq W_NV=R_1V$, and, since
$\dim_\field V=\dim_\field W_0V$, we have $\dim_\field R_1V\geq N + \dim_\field V$.
\end{proof}

We can improve on the bounds given in \eqref{naivebnds} by using more of the information available to us.
The key depends on the following lemma.

\begin{lem}\label{lem:cokbounds}
Let $Z\subset\pr2$ be a fat point subscheme. Assume
$t\geq\alpha(I_Z)$ and $h^1(\pr2, \I_Z(t))=0$ and let $e=h^1(\pr2, \I_Z(t-1))$. Then 
$e = \diff h_Z(t) = h_Z(t) - h_Z(t-1)$ and
$\max(2e-t,0)\leq \nu_{t+1}(Z)\leq e$.
\end{lem}

\begin{proof}
Let $L$ be a general line.
Taking cohomology of the short exact sequence
$0 \to \I_Z(t-1) \to \I_Z(t) \to \O_L(t) \to 0$
gives $h^0(\I_Z(t-1)) - h^0(\I_Z(t)) + (t+1) - e = 0$, since $h^1(\I_Z(t)) = 0$.
Rewriting in terms of $h_Z$ gives $e = \diff h_Z(t)$.

To bound $\nu_{t+1}(Z)=\dim\cok(\mu_{t}(Z))$, we first restrict to 
$L$ (defined, say, by a linear form $F$).
Twisting $0\to \O_{\pr2}(-L)\to \O_{\pr2}\to \O_L\to 0$ by $\O_{\pr2}(L)\cong \O_{\pr2}(1)$,
taking global sections 
and tensoring by $(I_Z)_t$ gives the short exact sequence
$0\to (I_Z)_t\to (I_Z)_t\otimes R_1\to (I_Z)_t\otimes H^0(L,\O_L(1))\to 0$.

We also have a second short exact sequence.
Since $L$ is general, we have $Z:L=Z$. 
Twisting \eqref{eqn:residual-ses} by $t+1$ and taking global sections
gives $0\to (I_Z)_t\to (I_Z)_{t+1}\to H^0(L, \O_L(t+1))\to 0$, which is a short exact sequence
since $h^1(\pr2, \I_Z(t))=0$.

The first short exact sequence maps to the second, giving the commutative diagram:

{
\[\begin{matrix}
0 & \to & (I_Z)_t & \to & (I_Z)_t\otimes R_1
& \to & (I_Z)_t\otimes H^0(L, \O_L(1)) & \to & 0 \\
  &     & \downarrow \mu_1 &     & \downarrow \mu_2 & & \downarrow \mu_3 & &  \\
0 & \to & (I_Z)_t & \to & (I_Z)_{t+1} 
& \to & H^0(L, \O_L(t+1)) & \to & 0 \\
\end{matrix}\]
}

Since $\mu_1$ is an isomorphism, we have $\cok\mu_2\cong\cok\mu_3$ by the snake lemma.
Twisting \eqref{eqn:residual-ses} by $\O_{\pr2}(t)$ and taking global sections gives
the four term exact sequence
$0\to (I_Z)_{t-1}\to (I_Z)_t\to H^0(L, \O_L(t))\to H^1(\pr2, \I_Z(t-1))\to 0$.
Let $V$ be the image of $(I_Z)_t$ in $H^0(L, \O_L(t))$.
Then $(I_Z)_t\otimes H^0(L, \O_L(1))$ and $V\otimes  H^0(L, \O_L(1))$
have the same image in $H^0(L, \O_L(t+1))$. By Lemma \ref{lem:genimfact},
the dimension of the image of $V\otimes  H^0(L, \O_L(1))$ is at least $\dim V+1$,
hence $\dim\cok\mu_2=\dim\cok\mu_3\leq h^0(L, \O_L(1))-(\dim V+1)=(t+2) - (\dim V+1)$.
Using the four term exact sequence, $\dim V=t+1-e$, so we have $\nu_{t+1}(Z)\leq e$.

Of course, the maximum possible dimension of the image of $\mu_3$
is obtained by assuming it has maximal rank. This gives $\dim\Image\mu_3\leq \min(2\dim V, t+2)$,
hence $\nu_{t+1}(Z)\geq (t+2)-\min(2\dim V, t+2)=\max(2e-t,0)$.
\end{proof} 

\begin{notation}\label{bettinotation}
Fixing notation, 
let $\dvec=(d_1,\ldots,d_{n+1})$ be a \GMS\ reduction vector for a
fat point subscheme $A$ with respect to a totally reducing sequence of lines
$L_1,\ldots,L_{n+1}$, and let the corresponding residual
schemes be $A=A_0,A_1,\ldots,A_n,A_{n+1}=\varnothing$.
For $\alpha(I_A) \leq t < \reg(I_A)$,
let $j = j(\dvec, t)$ be the greatest value in the interval $0 \leq j \leq n+1$
such that $h_{I_{A_i}}(t-i)=h_{I_A}(t)$ for $1\leq i\leq j$.
Note that each $A_i$ has a \GMS\ reduction vector given by a truncation of $\dvec$,
so each $h_{I_{A_i}}$ is determined by $\dvec$; thus $j$ is indeed determined by $\dvec$ and $t$.
Note also that $j \leq t$ (either $j \leq n+1 \leq t$; or if $t < n+1$ and $j > t$, then $A_t \neq \varnothing$,
so $0 = h_{I_{A_t}}(0) = h_{I_A}(t) > 0$, since $t \geq \alpha_A$).
\end{notation}

\begin{lem}\label{bettilemma}
With notation as in \ref{bettinotation}, we have:
\begin{enumerate}
\item $F = F_1 \dotsb F_j$ is a common divisor for $(I_A)_t$, where $F_i$ is a linear form defining $L_i$,
so that
\[
  \nu_{t+1}(A) = h_{I_A}(t+1) - h_{I_{A_j}}(t-j+1) + \nu_{t-j+1}(A_j) .
\]
\item The fat point scheme $A_j$ satisfies $h^1(\pr2, \I_{A_j}(t-j)) = 0$.
\end{enumerate}
\end{lem}

\begin{proof}Multiplication by $F$ gives an injective map $(I_{A_j})_{t-j} \to (I_A)_t$.
By assumption these have the same dimension, so also multiplication by $F$ is surjective, meaning $F$ is a common divisor for $(I_A)_t$.
Since multiplication by linear forms commutes with multiplication by $F$,
the image $\Image(\mu_t(A))$ of $\mu_t(A)$ in $(I_A)_{t+1}$ is $F \cdot \Image(\mu_{t-j+1}(A_j))$,
which has dimension $h_{A_j}(t-j+1) - \nu_{t-j+1}(A_j)$.
This gives $\dim \Image(\mu_t(A))$, so $\nu_{t+1}(A) = \codim \Image(\mu_t(A)) = h_A(t+1) - (h_{A_j}(t-j+1) - \nu_{t-j+1}(A_j))$.

If $j = n+1$ then $\I_{A_j} = \O_{\pr 2}$ and the vanishing of $h^1$ is automatic.
Otherwise, we claim $h^1(\pr2, \I_{A_i}(t-i)) = 0$ for $j \leq i \leq n+1$.
This is automatic for $i=n+1$.
In the short exact sequence
\[
  0 \to \I_{A_{i+1}}(t-i-1) \to \I_{A_i}(t-i) \to \O_{L_i}(t-i-d_i) \to 0 ,
\]
we may by downward induction assume $h^1(\pr2, \I_{A_{i+1}}(t-i-1)) = 0$,
and we have $t-i-d_i \geq t-j-d_j-1$ since $\dvec$ is \GMS.
From the short exact sequence
\[
  0 \to \I_{A_{j+1}}(t-(j+1)) \to \I_{A_j}(t-j) \to \O_{L_j}(t-j-d_j) \to 0
\]
the maximality of $j$ means $h_{I_{A_{j+1}}}(t-j-1) < h_{I_{A_j}}(t-j)$, whence $h^0(L_j,\O_{L_j}(t-j-d_j)) > 0$.
Thus $t-i-d_i \geq t-j-d_j-1 \geq -1$, and $h^1(\O_{L_i}(t-i-d_i)) = 0$.
This completes the downward induction.
\end{proof}

We can now give bounds improving on \eqref{naivebnds}.

\begin{thm}\label{betterbnds}
With notation as in \ref{bettinotation},
let $e = h^1(\pr2, \I_{A_j}(t-j-1)) = \diff h_{A_j}(t-j)$.
Then
\begin{equation*}
\begin{split}
\nu_{t+1}(A)\geq & (h_{I_A}(t+1)-h_{I_{A_j}}(t-j+1))+\max(2e-t+j,0),\\
\nu_{t+1}(A)\leq & (h_{I_A}(t+1)-h_{I_{A_j}}(t-j+1))+e .
\end{split}
\end{equation*}
\end{thm}
Note that the bounds for $\nu_{t+1}(A)$ are determined by $\dvec$ and $t$
(via truncations of $\dvec$ and $j = j(\dvec,t)$).

\begin{proof}
Applying Lemma \ref{lem:cokbounds}
with $Z=A_j$, $t-j$ in place of $t$ and $e=h^1(\pr2,\I_{A_j}(t-j-1))$
to the result of Lemma \ref{bettilemma}
gives the result.
\end{proof}

Since $t\geq j$, we see that $2e-t+j\leq 2e$, and hence
that the bounds in Theorem \ref{betterbnds} agree if $e=0$, and thus in that case we
even obtain an exact determination of $\nu_{t+1}(A)$. The following theorem gives a criterion for
the bounds on $\nu_t(A)$ to coincide for all $t$. A precise determination of exactly those $\dvec$ 
for which the bounds coincide for all $t$ can also be given; see Remark \ref{exactBettis}.

\begin{thm}\label{thm:exactbettis}
Let $A\subset\pr2$ be a fat point subscheme with full reduction vector $\dvec=(d_1,\ldots,d_{n+1})$.
If $\dvec$ is positive and strictly decreasing, then the bounds on $\nu_A(t+1)$ given in 
Theorem \ref{betterbnds} agree for all $\alpha(I_A)\le t< \reg(I_A)$, and thus $\dvec$ determines 
the graded Betti numbers for $I_A$.
\end{thm}

\begin{proof} By Proposition \ref{GMSprop}, $\dvec$ is \GMS. Thus by 
Theorem \ref{GMSthm}, Equation \eqref{tag3} and Definition \ref{thefunctions1}, we have
$h_A=f_\dvec$. By Remark \ref{thefunctions4}, $f_\dvec(t)$ is the number of dots of
$S_\dvec$ on or to the left of the line $x+y=t$. But $\alpha(I_A)$ is the least $t$ such that
$f_\dvec(t)=h_A(t)<\binom{t+2}{2}$; i.e., it is the least $t$ such that there is
an integer lattice point $(i,j)\not\in S_\dvec$ with $i\ge 0$, $j\ge0$ and $i+j\le t$.
Since the entries of $\dvec$ are strictly decreasing, this is precisely the
number of rows of $S_\dvec$; i.e., $\alpha(I_A)=n+1$.
And the regularity $\reg(I_A)$ is $t+1$ for the least degree $t$ such that
$h_A(t)=\deg(A)$, which is just the least $t$ such that every dot of $S_\dvec$ is
on or to the left of the line $x+y=t$. Since $\dvec$ is strictly decreasing,
this $t$ is $d_1-1$; i.e., $\reg(I_A)=d_1$. 

Now suppose that $\alpha(I_A)\le t<\reg(I_A)$, and let $L_1,\ldots,L_{n+1}$ be the sequence
of lines giving rise to the reduction vector $\dvec$. The divisor $L_1+\cdots+L_j$ is in 
the base locus of $(I_A)_t$ if $h_{I_A}(t)=h_{I_{A_i}}(t-i)$ for $0< i\le j$. 
Assume we pick the largest $j$ such that this holds. This is equivalent to saying that
the line $x+y=t$ contains dots of $S_\dvec$ for each of the rows of $S_\dvec$
corresponding to $d_1,\ldots,d_j$ but for none of the other rows.
In particular, $t-j\ge \reg(I_{A_j})$ and hence $e=h^1(\pr2,\I_{A_j}(t-j-1))=0$,
so, as pointed out above following Theorem \ref{betterbnds}, our bounds coincide.
\end{proof}

\begin{example}
Let $Z=3(p_1+\cdots+p_{10})$ where the points $p_i$ are the points of pairwise
intersections of five general lines $L_1,\ldots,L_5$ (as shown in Figure~\ref{star5fig}).
Using the sequence of lines $L_1,\ldots,L_5,L_1,\ldots,L_4$ where $L_i$ is as shown
in the figure, we obtain a full reduction vector $\dvec=(12,11,10,9,8,4,3,2,1)$ for $Z$.
This is \GMS\ so we obtain an exact determination $h_Z=(1,3,6,10,15,21,28,36,45,50,55,60,60,\ldots)$
for the Hilbert function of $Z$. From this we see that $\alpha(I_Z)=9$ and $\reg(I_Z)=12$,
hence $\nu_t(Z)=0$ for $t<9$ and for $t>12$, and $\nu_9(Z)=\binom{9+2}{2}-h_Z(9)=5$.
For $t=9$, the bounds given in \eqref{naivebnds} 
are $0=\max(0,11-3\cdot5)\leq \nu_{10}(Z)\leq11-(2+5)=4$. In this case $e=0$,
so Theorem \ref{betterbnds} gives the exact value 
$\nu_{10}(Z)=h_{I_Z}(t+1)-h_{I_{Z_j}}(t-j+1)=0$.

Similarly, $t=10$ gives $\nu_{11}(Z)=0$; now consider $t=11$. Then \eqref{naivebnds}  gives the bounds
$0=\max(0,31-3\cdot18)\leq \nu_{12}(Z)\leq31-(2+18)=11$, but 
from Theorem \ref{betterbnds} again $e=0$ so we have the exact value, 
$\nu_{12}(Z)=h_{I_Z}(t+1)-h_{I_{Z_j}}(t-j+1)=31-26=5$.
\end{example}

\begin{remark}\label{exactBettis}
We have indicated above how to find bounds in terms of a full reduction vector for the graded Betti numbers for 
fat point subschemes of $\pr2$ when the reduction vector is \GMS. 
In the preceding example we saw that we obtained exact values for $\nu_t$ for every $t$,
as we expect from Theorem  \ref{thm:exactbettis} since the reduction vector is positive and
strictly decreasing. It is possible to show that our approach gives
an exact determination for $\nu_t$ for all $t$ for a positive \GMS\ reduction vector $\dvec$
whenever (and only when) the reduction vector
has one of the following forms (see \cite[Proposition 5.1.8]{refCHT}):
\begin{enumerate}
\item $\dvec$ is strictly decreasing;
\item $\dvec=(m,m, m-1,\dots,2,1)$ where $m \geq 1$;
\item $\dvec = (d_1,\dotsc,d_k,m,m,m-1,\dotsc,2,1)$ where $(d_1,\dotsc,d_k)$ is strictly decreasing and $d_k \geq m+2$.
\end{enumerate}
Because the proof is somewhat long but does not involve any 
new ideas, we have omitted the proof here.
\end{remark}

\begin{remark}\label{GMS-CHTcomparison} 
Our criterion for when a reduction vector $\dvec$ determines the Betti numbers of a fat point scheme
is more restrictive than an analogous criterion given in \cite{refGMS}.
This is because the criterion given in \cite{refGMS} is only for the situation of linear configurations
(see section \textsection\ref{linconfig} for the definition), whereas we allow arbitrary fat point schemes.
We give an example to illustrate this.
Let $\dvec=(2,2)$.
Then $\dvec$ satisfies the criterion given in \cite[Theorem 4.5]{refGMS} 
for the graded Betti numbers of $I_A$ for a linear configuration $A$ of type $\dvec$
to be uniquely determined.
Indeed, $A$ being a linear configuration of type $(2,2)$ means simply that $A$ is a reduced
set of $4$ points in general position.
In particular $\dvec=(2,2)$ is a reduction vector for $A$, obtained by reducing along two lines
through pairs of the points.

Our methods give bounds for the graded Betti numbers of $I_B$
for any fat point scheme $B$ having reduction vector $\dvec=(2,2)$, since $(2,2)$ is \GMS.
However, as indicated by Remark \ref{exactBettis}, these bounds do not coincide, 
so they do not determine the graded Betti numbers of $I_B$.
Indeed without assuming $B$ is a linear configuration of type $(2,2)$ the graded Betti numbers
are not uniquely determined---even assuming $B$ is reduced is not sufficient.
For instance, let $B$ be the reduced union of $4$ points with exactly $3$ collinear.
Let $L_2$ be any line determined by the fourth (non-collinear) point and any one
of the collinear points, and let $L_1$ be the line through the other two points.
We then obtain the same reduction vector $\dvec=(2,2)$ by taking
our reducing sequence of lines to be $L_2, L_1$, in that order.
But since $B$ contains $3$ collinear points, the graded Betti numbers of $I_B$ are different
from those of $I_A$.
(Note that $L_1, L_2$ gives $(3,1)$ as a reduction vector for $B$, determining the graded Betti numbers of $I_B$.
Indeed this $B$ is a linear configuration of type $(3,1)$.)

In conclusion, since the bounds we obtain must
be broad enough to accommodate subschemes with different graded Betti numbers
but still having the same $\dvec$, there are situations where our bounds will not coincide
(such as the case described here of four general points) for which \cite{refGMS}, 
as a result of restricting to the case of linear configurations,
can give an exact determination of the Betti numbers.
\end{remark}

\begin{proof}[Proof of Theorem \ref{mainthmintro}]
This theorem is an immediate consequence of  Corollary \ref{idealbounds} and Theorems \ref{GMSthm} and \ref{thm:exactbettis},
in view of Proposition \ref{GMSprop}, Equations \eqref{tag1} and \eqref{tag3}, and Definition \ref{thefunctions1}.
\end{proof}

\section{Applications and Examples}\label{examplessection}

In this section we present a collection of examples.

\subsection{Linear and line count configurations}\label{linconfig}

We first consider some examples which generalize one of the main results of \cite{refGMS}.  The following notation is needed for these examples.

\begin{defn}
Given an integer vector $\vvec=(v_1,\ldots,v_n)$,
we define the \defining{permuting operator} $\pi(\vvec)$ to be the vector
whose entries are the entries of $\vvec$ permuted to be non-increasing.
\end{defn}

\begin{defn}\label{staropdef}
Given positive integer vectors $\avec=(a_1,\ldots,a_n)$ and $\mvec=(m_1,\ldots,m_n)$, let
\[ \avec \circ \mvec = (a_1m_1,(a_1-1)m_1,\ldots,2m_1, m_1,\hspace{1ex} a_2m_2,(a_2-1)m_2,\ldots,m_2,\hspace{1ex} \ldots, 
\hspace{1ex} a_nm_n,(a_n-1)m_n,\ldots,m_n). \]
The \defining{star operator} is defined by $\avec * \mvec =\pi(\avec \circ \mvec)$.
Note that the number of entries in $\avec \circ \mvec$ or $\avec * \mvec$ is $\sum_{i=1}^n a_i$.
\end{defn}

\begin{example}
We have
\[(2,3)*(2,3) = \pi((4,2,\hspace{1ex} 9,6,3))=(9,6,4,3,2)\] 
and 
\[(3,2)*(2,3) = \pi((6,4,2,\hspace{1ex} 6,3))=(6,6,4,3,2).\]
\end{example}

\begin{example}
Note that $(3,1)*(2,3) = (2,2)*(2,3) = (6,4,3,2)$, so we cannot recover $\avec$ from $\avec * \mvec$ even if we know $\mvec$.
Similarly, $(1,2)*(2,4) = (1,2)*(8,2)=(8,4,2)$, so we cannot recover $\mvec$ from $\avec * \mvec$ even if we know $\avec$.
\end{example}

One of the motivating goals of the paper \cite{refGMS} was the determination (in characteristic 0) of the
Hilbert function for $2Z$, whenever $Z \subset \pr2$ is a \defining{linear configuration of type $(m_1, \ldots, m_s)$};
i.e., whenever $Z=Z_1+\dots+Z_s$ for some set of distinct lines $L_1,\dots,L_s$ 
and finite reduced schemes $Z_i\subset L_i$, such that no point of $Z$
occurs where two of the $L_i$ meet,
and the integers $m_i=\card(Z_i)$, $i=1,\dots,s$ are distinct. 
After re-indexing, we may assume that $m_1>\dots>m_s>0$.  
In our terminology, \cite{refGMS} showed that $h_{2Z} = \sum \diag(\avec * \mvec)$
where $\mvec = (m_1,\dots,m_s)$ and $\avec = (2,\dots,2)$, if $\avec * \mvec$ is \GMS.

The following corollary is a generalization of the above result of \cite{refGMS} in which
we give an exact determination of the Hilbert function in a broader range of cases.
If we allow the integers $m_i$ to be equal in the definition of a linear configuration,
then we obtain what we call a \defining{line count configuration}.  

\begin{thm}\label{GMSgeneralization} 
Given positive integers $a_i$ and $m_i$, let $\avec = (a_1, \dots, a_s)$ and $\mvec=(m_1, \dots, m_s)$.
Let $Z_1 + \cdots + Z_s$ be a line count configuration of type $\mvec$ corresponding to distinct lines
$L_1, \ldots, L_s$ and let $Z(\avec,\mvec) = a_1Z_1 + \cdots + a_sZ_s$.
Then $\avec * \mvec$ is a full reduction vector for $Z(\avec,\mvec)$, hence
$f_{\avec * \mvec}\leq h_{Z(\avec, \mvec)}\leq F_{\avec * \mvec}$.
If $\avec * \mvec$ is \GMS, then $h_{Z(\avec,\mvec)}=f_{\avec * \mvec}=F_{\avec * \mvec}=\sum \diag(\avec * \mvec)$.
\end{thm}

\begin{proof} Consider the sequence of lines ${\bf L}=\{L_1^{a_1},\dots,L_s^{a_s}\}$, where 
the notation $L_i^{a_i}$ denotes the sequence $L_i,\dots,L_i$ with $a_i$ terms.
It is clear that ${\bf L}$ is a totally reducing sequence of lines for $Z(\avec,\mvec)$ with reduction vector
$\avec \circ \mvec$. 
By definition, $\avec*\mvec$ is obtained from $\avec \circ \mvec$ by appropriately
permuting the entries of $\avec \circ \mvec$. By applying the same permutation to
the entries of ${\bf L}$, we obtain another totally reducing sequence of lines for $Z(\avec,\mvec)$ but with reduction vector
$\avec * \mvec$. Thus $f_{\avec * \mvec}\leq h_{Z(\avec,\mvec)}\leq F_{\avec * \mvec}$, and
if $\avec * \mvec$ is \GMS, then $f_{\avec * \mvec}=F_{\avec * \mvec}$, and hence  
$h_{Z(\avec,\mvec)}=f_{\avec * \mvec}=\sum \diag(\avec * \mvec)$.
\end{proof}

Our methods determine the Hilbert functions of fat point schemes whose support are line count configurations consisting of one and two 
lines, as the next example shows.

\begin{example}\label{lines}
Let $Z$ consist of $m$ points on a line $\ell$ in $\pr2$ and define the fat point scheme $A := aZ$.  Then the sequence of lines 
$\ell^a$  (where again $\ell^a$ denotes the sequence $\ell, \dots, \ell$ in which $\ell$ is repeated $a$ times)
is a totally reducing sequence for $A$ with \GMS\ reduction vector $(a)*(m)= (am,\dots, 2m,m)$.  
Thus $h_A = f_{(a)*(m)} = F_{(a)*(m)} = \sum \diag((a) * (m))$.

Similarly, suppose $Z = Z_1 + Z_2$ where $Z_1$ is a set of $m_1$ points
on a line $\ell_1$ and $Z_2$ is a set of $m_2$ points on a line $\ell_2$ such that
no point of $Z$ lies at the intersection of $\ell_1$ and $\ell_2$.
Assume that $m_1 \geq m_2 > 2$.
Consider the fat point scheme $A = a_1Z_1 + a_2Z_2 \subset \pr2$.
The sequence of lines $\ell_1^{a_1}, \ell_2^{a_2}$ is a totally reducing sequence
of lines for $A$ with reduction vector $\avec \circ \mvec$, where $\avec = (a_1, a_2)$
and $\mvec = (m_1, m_2)$.
Again, after a permutation of the lines, we obtain a totally reducing sequence of lines
with reduction vector $\avec*\mvec$.
But $\avec * \mvec$ is \GMS\ by Proposition \ref{GMSprop} (since at most two consecutive entries 
of $\avec * \mvec$ can be equal and otherwise they differ by at least 2).
Thus, again, $h_A = f_{\avec * \mvec} = F_{\avec * \mvec} = \sum \diag(\avec * \mvec)$.

For fat point schemes supported on line count configurations involving three or more lines,
these arguments do not always determine the Hilbert functions.
For example, let $Z = Z_1 + Z_2 + Z_3$, where 
each
$Z_i$ is a set of $m \geq 3$ points
on a line $\ell_i$ such that no point of $Z$ lies in the intersection of the lines.
Consider the fat point scheme $A = 2Z_1 + 2Z_2 + 2Z_3\subset \pr2$.
The sequence of lines $\ell_1,\ell_2 ,\ell_3,\ell_1,\ell_2,\ell_3$
is a totally reducing sequence of lines for $A$ with reduction vector $\avec * \mvec$
where $\avec = (2,2,2)$ and $\mvec = (m, m, m)$.
However, $\dvec=\avec * \mvec$ is not \GMS, so $f_\dvec(t)<F_\dvec(t)$
for at least one value of $t$.
\end{example}

\subsection{Intersections of lines}

Two things (at least) make linear configurations interesting.
One is that it is easy to write down their Hilbert functions:
a (reduced) linear configuration $Z$ of type ${\bf c}=(c_1,\ldots,c_s)$
has Hilbert function $\sum\diag({\bf c})$. Another is that according to
results of Macaulay \cite{refMa} and Geramita-Harima-Shin \cite{refGHS}, every Hilbert function
for 
0-dimensional
subschemes of $\pr2$ already occurs for
linear configurations.

However, when considering Hilbert functions of schemes $aZ$
where $Z$ is a reduced 0-dimensional subscheme of $\pr2$,
linear configurations do not behave well. In fact,
if $Z$ is a linear configuration of type ${\bf c}=(c_1,\ldots,c_s)$
for $s>2$ which we assume for simplicity has $c_1> c_2> \cdots> c_s$,
and if $a\ge c_{s-2}c_{s-1}$, then $(a,\ldots,a)* {\bf c}$ will have an entry
equal to $c_{s-2}c_{s-1}c_s$ which occurs at least three times, and hence
is not \GMS, and so our upper and lower bounds on
the Hilbert function of $aZ$ will not coincide.

Linear configurations arise by taking finite sets of points on a union of lines
$L_1,\ldots,L_s$ but avoiding the points of intersections of the lines.
It turns out the opposite extreme, of taking $Z$ to be the set of all of the points
where two or more lines $L_i\ne L_j$ intersect, behaves better in this respect.
In fact, if in this situation the number of lines $L_i$ through each point of $Z$
is the same (say $n$ lines go through each point of $Z$), then our results determine the
Hilbert function of $mnZ$ and the graded Betti numbers of $I_{mnZ}$ for {\it every} $m\ge 1$.
In particular, 
Example \ref{sArbLinesButn=2} handles the case $n=2$, 
Example \ref{PrimeCharEx} handles the case $n=q$ where $q$ is a power of a prime, and
Example \ref{HesseEx} handles the case $n=3$. 
If $n=2$, then we even obtain the Hilbert function of $mZ$ and the graded Betti
numbers of $I_{mZ}$ for {\it every} $m\ge 1$.

Here is the general set up. Let $L_1,\ldots,L_s$ be distinct lines in $\pr2$, and let
$e_1,\ldots,e_s$ be positive integers (these need not be in non-decreasing order).
Consider the divisor $D=\sum_ie_iL_i$.  Let
$S_D$ be the finite set of all points $p$ which lie simultaneously on
two or more of the lines $L_i$, and define $Z(D)$ to be the fat point scheme
$\sum_{p\in S_D}m_pp$,
where $m_p = \sum_{p \in L_i} e_i$
is the number of lines $L_i$ which pass through $p$, counting with multiplicity.
In case $e_1=\cdots=e_s=1$ (in which case $D$ is reduced), we also define
$Z'(D)=\sum_{p\in S_D}(m_p-1)p$.
In Example~\ref{example1: reduction of intersections of s arbitrary lines}
we consider $Z'(D)$ for $D$ reduced.
In Example~\ref{sArbLines} we consider $mZ(D)$ for any $m$ and any $D$.
In Example \ref{sArbLinesButn=2}, under the special condition
that no three lines of $D$ have a common point, we consider $mZ'(D)$ for any $m$ and for $D$ reduced.

In each case we describe the reduction vectors,
which turn out to be not only \GMS, but even strictly decreasing.
Hence the Hilbert functions are completely determined by the reduction vector and,
as pointed out in Remark \ref{exactBettis}, so are
the graded Betti numbers.

For comparison, Ardila and Postnikov \cite{refAP} obtain the Hilbert function 
both for $Z'(D)$ and for $Z(D)$, not only for $\pr2$ but for hyperplane arrangements
in any dimension, over the complex numbers.
They express their solution in terms of the Hilbert series.
While on the one hand our results are only for subschemes of $\pr2$, on the other hand
for us the characteristic is arbitrary, and we obtain
the Hilbert functions and Betti numbers, not only for $Z'(D)$ and for $Z(D)$ but also
as mentioned above for $mZ(D)$ (in Example \ref{sArbLines}) and 
for $mZ'(D)$ (in Example \ref{sArbLinesButn=2}).

\begin{example}\label{example1: reduction of intersections of s arbitrary lines}
Consider $Z'_{s-1}=Z'(D)$, where $D=L_1+\cdots+L_s$ for distinct lines $L_i$, so $D$ is reduced.
Note that the scheme theoretic intersection
$Z'_{s-1}\cap L_i$ of any of the lines $L_i$
with $Z'(D)$ has degree exactly $\sum_{p\in L_i\cap S_D} (m_p-1)$, but this sum is just
the number of lines other than $L_i$; i.e., it is $s-1$.
The residual $Z'_{s-2}$ of $Z'_{s-1}$ with respect to
$L_i$ is $Z'(D-L_i)$.
Then for any line $L_j$ with $j\ne i$,
the intersection $Z'_{s-2}\cap L_j$ has degree $s-2$.
Thus taking residuals with respect to the sequence of lines $L_2,\ldots,L_s$
(or indeed any sequence of $s-1$ of the $s$ lines)
results in a
full reduction vector
$\dvec=(s-1,s-2,\ldots,3,2,1)$.
\end{example}

\begin{example}\label{sArbLines}
Now consider $Z(D)$ for $D=e_1L_1+\cdots+e_sL_s$ for distinct lines $L_i$ and $e_i\ge1$, so $D$ need not be reduced.
Whereas in Example \ref{example1: reduction of intersections of s arbitrary lines}
we could totally reduce $Z'(D)$ and obtain a strictly decreasing (and hence \GMS) reduction vector
without regard to the order in which we chose the lines (as long as we
chose $s-1$ different lines), now the order will matter, but 
we can nevertheless obtain a strictly decreasing reduction vector not only for $Z(D)$
but for $mZ(D)$ for any $m>0$.

In general, the reduction vector corresponding to an arbitrary choice of
a sequence of lines need not be \GMS, or even non-increasing; we now show how
the sequence of lines can be chosen so that in fact we do get a strictly decreasing reduction vector,
by using the following ``greedy algorithm'':
at any given step, say the residual fat point scheme is $Z_k=\sum_{p\in S_D}a_pp$.
Choose any one of the lines $L_i$ maximizing the degree of $Z_k\cap L_i$
(i.e., maximizing the sum $\sum_{p\in L_i\cap S_D} a_p$),
among lines which have not been chosen $me_i$ times so far.
This is totally reducing, as for each point $p \in S_D$, each line $L_i$ through $p$
occurs $m e_i$ times in the reducing sequence, giving a total
of $\sum m e_i = m m_p$ lines through $p$ (where the sum is over all $i$ such that $p\in L_i$).

We now explain why the resulting reduction vector in fact contains no repeated values.
If the sum of the multiplicities at the $k$th step along some line, say $L_i$,
is the maximum, then in the next step, the sum of the multiplicities along any line, say $L_j$,
is strictly less, because the multiplicities along line $L_i$ were reduced
and for any $j\ne i$ the multiplicity of the point $L_i \cap L_j$ has been reduced by $1$,
as long as that multiplicity was not already 0.
But for the multiplicity to have already been 0, all of
the lines through that point must have been chosen already; in particular,
$L_j$ must have been chosen $me_j$ times already, and
so we do not allow ourselves to choose it again.
\end{example}

\begin{figure}[h]
\setlength{\unitlength}{12mm}   
\begin{center}
\begin{picture}(10,3)(-3.5,-.5)
\linethickness{.25mm}
\put(0.1,0.1){\line(1,2){1.5}}
\put(0.875,3.375){\line(1,-3){1.35}}
\put(-1.25,2.5){\line(5,-4){3.7}}
\put(-.3,.2){\line(5,3){3.2}}
\put(-1,2.1){\line(6,-1){4}}
\put(0.344,0.583){\circle*{.2}}
\put(2.045,-.14){\circle*{.2}}
\put(2.03,1.6){\circle*{.2}}
\put(1.225,2.328){\circle*{.2}}
\put(-.68,2.048){\circle*{.2}}
\put(0.809,0.867){\circle*{.2}}
\put(1.56,1.31){\circle*{.2}}
\put(1.436,1.7){\circle*{.2}}
\put(.945,1.788){\circle*{.2}}
\put(.57,1.04){\circle*{.2}}
\put(1.5,2.6){$L_{1}$}
\put(2.5,1.21){$L_{2}$}
\put(1.75,-.6){$L_{3}$}
\put(-.43,.47){$L_{4}$}
\put(-1,2.33){$L_{5}$}
\end{picture}
\end{center}
\caption[The configuration $Z(s)$ with $s=5$.]{The configuration $Z(L_1 + L_2 + L_3 + L_4 + L_5)$.}\label{star5fig}
\end{figure}
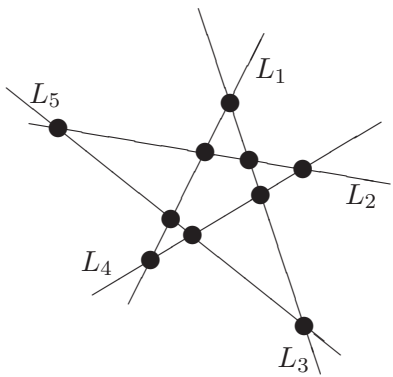

\begin{example}\label{sArbLinesButn=2}
Here we consider so-called {\it star configurations}; see Figure~\ref{star5fig}.
These are of the form $Z'(D)$ in the special case that $D=L_1+\cdots+L_s$ is 
a sum of distinct lines (so $D$ is reduced) under the assumption that
at any point $p=L_i\cap L_j$, $i\ne j$, where two lines meet, these are the only
two lines which contain $p$. Such configurations arise in
an important way in both \cite{refGMS} and \cite{refBH} as having extremal behavior.
In \cite{refGMS},
$h_{mZ'(D)}$ is determined for $m=2$, while \cite{refBH} determines $\alpha(I_{mZ'(D)})$ for every $m$.
Our results give both $h_{mZ'}$ and the graded Betti numbers for $I_{mZ'}$ since
the reduction vector is positive and strictly decreasing, for all $m>0$, as we now show.
Note that $Z'=Z'(D)$ is just the reduced union of the $\binom{s}{2}$ points
and $Z(D)=2Z'$.
Thus we get the Hilbert function and graded Betti numbers for $mZ'$ when $m$ is even by
Example \ref{sArbLines}, in which case the reduction vector is 
\[
\begin{split}
\big(m(s-1), m(s-1)-1, \ldots, m(s-1)-(s-1),\\
(m-2)(s-1), (m-2)(s-1)-1, \ldots, (m-2)(s-1)-(s-1),\\
\ldots,\\
2(s-1),2(s-1)-1,\ldots,2(s-1)-(s-1)\big),
\end{split}
\]
where we have used a multi-line display to make the pattern easier to discern.
We also get them when $m$ is odd, as follows.
Say $m=2a+1$ is odd. Then $mZ'=aZ(D)+Z'(D)$. First as we just showed we can choose a reducing
sequence of lines for $aZ(D)$, giving a strictly decreasing reduction vector
$\dvec=(d_1,\ldots,d_u)$ as displayed in the even case, 
then choose a reducing sequence of lines for $Z'(D)$,
giving the reduction vector $\dvec'=(s-1,\dots,2,1)$.
Putting these sequences together reduces $mZ'$ first to $Z'(D)$,
and then it reduces $Z'(D)$ to the empty scheme, and the reduction
vector is $(d_1+s-1,\dots,d_u+s-1, s-1,\ldots,2,1)$,
which is clearly positive and strictly decreasing.
\end{example}

\vskip2\baselineskip

\addtocounter{thm}{1}

\begin{example}\label{PrimeCharEx}
Here we show that in positive characteristic there can be examples
where the number of lines through each point is larger than 2.
Let $Z$ be the reduced union of the points
of $\PPP^2_{\field}$, where $q$ is a power of a prime $p$ and 
$\field={\bf F}_q$ is the finite field of $q$ elements.
For the set of lines, take all of the lines in $\PPP^2_{{\bf F}_q}$. Then
there are $q^2+q+1$ points, $q^2+q+1$ lines, $q+1$ lines go through
every point, and $q+1$ points lie on every line. If $D$ is any divisor
with support on all of the lines, then by Example \ref{sArbLines}
we obtain $h_{mZ(D)}$ and the graded Betti numbers for $I_{mZ(D)}$,
for all $m>0$.
In particular, we obtain $h_{m(q+1)Z}$ and the graded Betti
numbers for $I_{m(q+1)Z}$, for all $m>0$.
\end{example}

\begin{example}\label{HesseEx}
Finally we consider a classical example over the complex numbers where the number of lines
through each point is 3.
In particular, consider the projective dual of the Hesse configuration
where $\field$ is the complex numbers
(see \cite{refTot} for another situation in which this has come up recently).
The Hesse configuration is the set of nine flex points for a general smooth plane cubic $C$
and the twelve lines through 
each pair of
flexes.
It has the property that any of the twelve lines through two flexes goes through a third,
and each of the nine flexes lies on four of the lines.
However, it is not true that every pair of the twelve lines meet at a flex,
so the Hesse configuration is not of the type we have considered in the previous examples;
i.e., it is not of the form $Z(D)$ for a divisor $D$ supported on the twelve lines.

Consider now the configuration projectively dual to the Hesse configuration,
consisting of twelve points and nine lines,
with each point lying on three of the lines and each line passing through four of the points.
Now each pair of the lines (corresponding to a pair of flex points on $C$)
does intersect at one of the twelve points (corresponding to a line through the two flex points on $C$).
For any divisor $D$ supported on these nine lines we get a reduction vector
for the fat point scheme $Z(D)$.

In particular, if $D$ is the reduced union of the nine lines and $Y$ is the reduced union
of the twelve points, then $Z(D) = 3Y$ and $Z'(D) = 2Y$.
By the previous examples we get strictly decreasing reduction vectors for
$m3Y$, $2Y$, and even $(m3+2)Y$, for all $m > 0$,
by an argument similar to that in Example \ref{sArbLinesButn=2}.
In fact, in this case one can also check ad hoc that $(3m+1)Y$ has a \GMS\
(but not strictly decreasing) reduction vector.
So we get the Hilbert function of $mY$ for all $m>0$,
and the graded Betti numbers for all $m$ congruent to $0$ or $2$ modulo $3$.
\end{example}


\begin{thebibliography}{GHM}

\bibitem[AP]{refAP} F.\ Ardila and A.\ Postnikov. {\it Combinatorics and Geometry of Power Ideals}, 
Trans. Amer. Math. Soc. 362 (2010), 4357--4384.

\bibitem[BH]{refBH} C.\ Bocci and B.\ Harbourne. {\it Comparing 
Powers and Symbolic Powers of Ideals}, 
Journal of Algebraic Geometry, {\bf 19} (2010) 399--417.

\bibitem[CM1]{refCM1} C.\ Ciliberto and R.\ Miranda. {\it
The Segre and Harbourne-Hirschowitz Conjectures},
in: Applications of algebraic geometry to coding theory, physics and
computation (Eilat 2001), NATO Sci. Ser. II Math. Phys. Chem., 36, Kluwer Acad. Publ.,
Dordrecht, (2001), 37 -- 51.

\bibitem[CM2]{refCM2} C.\ Ciliberto and R.\ Miranda. {\it
Linear Systems of Plane Curves with Base Points of
Equal Multiplicity}, Trans. Amer. Math. Soc. 352 (2000), 4037--4050.

\bibitem[CM3]{refCM3} C.\ Ciliberto and R.\ Miranda, {\it Nagata's Conjecture for a Square Number of Points}, 
Ricerche di Matematica, 55 (2006), 71--78. 

\bibitem[CHT]{refCHT} S.\ Cooper, B.\ Harbourne and Z.\ Teitler, 
{\it Using residuation and collinearity to bound Hilbert functions of fat points in the plane}, preprint, 	arXiv:0912.1915 version 1.

\bibitem[Du]{refDu} M.\ Dumnicki, {\it Reduction method for linear systems of 
plane curves with base fat points}, Ann. Polon. Math. 90.2 (2007), 131--143.

\bibitem[Ev1]{refEv} L.\ Evain, {\it La fonction de Hilbert de la r\'eunion de $4^h$ gros points g\'en\'eriques 
de $\pr2$ de m\^eme multiplicit\'e}, J.\ Algebraic Geom.\ 8 (1999), 787--796. 

\bibitem[Ev2]{refEv2} L.\ Evain, {\it Computing limit linear series with infinitesimal methods}, 
Ann.\ Inst.\ Fourier (Grenoble) 57 (2007), no. 6, 1947--1974. 

\bibitem[FHL]{refFHL} G.\ Fatabbi, B.\ Harbourne and A.\ Lorenzini. {\it Resolutions 
of Ideals of Fat Points with Support in a Hyperplane}, Proc.\ Amer.\ Math.\ Soc.\ 134 (2), 2006, 3475--3483.

\bibitem[FHH]{refFHH} S.\ Fitchett, B.\ Harbourne and S.\ Holay.
{\it Resolutions of Fat Point Ideals involving Eight General Points of $\pr2$}, 
J. Algebra 244 (2001), 684--705.

\bibitem[FL]{refFL} S.\ Franceschini and A.\ Lorenzini.  {\it Fat Points of ${\bf P}^n$ Whose Support is 
Contained in a Linear Proper Subspace},  J.\ Pure Appl.\ Alg.\ 160 (2-3), 169--182 (2001).

\bibitem[GHM] {refGHM} A.\ V.\ Geramita, B.\ Harbourne and J.\ Migliore. {\it 
Classifying Hilbert functions of fat point subschemes in $\pr2$},
Collect. Math. 60, 2 (2009), 159--192.

\bibitem[GHS]{refGHS} A.\ V.\ Geramita, T.\ Harima and Y.\ S.\ Shin.  {\it An Alternative 
to the Hilbert Function for the Ideal of a Finite Set of Points in ${\bf P}^n$}, Illinois J.\ Math.\ 45 (2001), no.\ 1, 1--23.

\bibitem[GMR] {refGMR} A.\ V.\ Geramita, P.\ Maroscia and L.\ G. Roberts.  {\it The Hilbert Function of a Reduced $\field$-Algebra},
J.\ London Math.\ Soc.\ (2), 28 (1983), 443--452.

\bibitem[GMS] {refGMS} A.\ V.\ Geramita, J.\ Migliore and L.\ Sabourin. {\it 
The First Infinitesimal Neighborhood of a Linear Configuration of Points in 
$\pr2$}, J. Algebra 298 (2), 2006, 563--611.

\bibitem[G]{refGi} A.\ Gimigliano. {\it On linear systems of plane curves}, Thesis, QueenÕs University, Kingston, 1987. 

\bibitem[H3]{refKrakow} B.\ Harbourne. {\it  Global aspects of the geometry of surfaces}, 
Ann. Univ. Paed. Cracov. Stud. Math. 9 (2010), 5--41
(\href{http://arxiv.org/abs/0907.4151}{arXiv:0907.4151}).

\bibitem[H4]{refVanc} B.\ Harbourne. {\it  The geometry of rational surfaces and 
Hilbert functions of points in the plane}, 
Canadian Mathematical Society Conference Proceedings 6, 95--111 (1986).

\bibitem[HR]{refHR} B.\ Harbourne and J.\ Ro\'e. {\it  Linear systems with 
multiple base points in $\pr2$}, Adv.\ Geom. 4 (2004), 41--59. 

\bibitem[HHF]{refHHF} B.\ Harbourne, S.\ Holay and S.\ Fitchett. 
{\it Resolutions of ideals of quasiuniform fat point subschemes of 
$\pr2$}, Trans. Amer. Math. Soc. 355 (2003), no. 2, 593--608.

\bibitem[Hi1]{refHi} A.\ Hirschowitz.
{\it La m\'ethode d'Horace pour l'interpolation \`a plusieurs
variables}, Manus. Math. 50 (1985), 337--388.

\bibitem[Hi2]{refHi2} A.\ Hirschowitz.
{\it Une conjecture pour la cohomologie des diviseurs sur les surfaces
rationelles g\'en\'eriques}, J. Reine Angew. Math., 397 (1989), 208--213.

\bibitem[I]{refI} A. Iarrobino. {\it Inverse system of a symbolic power III: thin algebras and fat points}, 
Compositio Math. 108, (1997), 319--356.

\bibitem[Ma]{refMa} F. S.\ Macaulay.  {\it Some Properties of Enumeration in 
the Theory of Modular Systems}, Proc. London Math. Soc.\ (2), 26 (1927), 531--555.

\bibitem[Mi]{refMi}  J.\ Migliore. {\it The Geometry of Hilbert Functions}, in: ``Syzygies and Hilbert 
functions,'' Lect. Notes in Pure and Appl. Math. 254 (ed. Irena Peeva), CRC Press (2007), 179--208.

\bibitem[N]{refNag} M.\ Nagata. {\it On rational surfaces, II}, Mem.\ Coll.\
Sci.\ Univ.\ Kyoto, Ser.\ A Math.\ 33 (1960), 271--293. 

\bibitem[R]{refR} J.\ Ro\'e, {\it Limit linear systems and applications}, 
arXiv:math/0602213, preprint, 2006.

\bibitem[Sc]{refHS} H.\ Schenck. {\it Geometry and Syzygies of Rational Surfaces Arising from Line Configurations in $\pr2$}, preprint, 2009.

\bibitem[St]{refS} R.\ Stanley.  {\it Hilbert Functions of Graded Algebras}, Adv. in Math. 28, 1978, 57 - 83.

\bibitem[To]{refTot}
B.\ Totaro, {\it The cone conjecture for {C}alabi-{Y}au pairs in dimension
  2}, Duke Math. J. \textbf{154} (2010), no.~2, 241--263.



\end{thebibliography}
\end{document}